\documentclass[preprint]{elsarticle}
\usepackage{amsmath,amssymb,color,float}
\usepackage[dvipsnames]{xcolor}
\usepackage{caption}
\usepackage{subcaption}

\numberwithin{equation}{section}
\numberwithin{figure}{section}
\numberwithin{table}{section}

\newtheorem{remark}{Remark}

\renewcommand{\Re}{\mathop{\rm Re}}

\newcommand{\bm}[1]{{\boldsymbol #1}}
\newcommand{\supp}{\mathop{\rm supp}}
\newcommand{\diag}{\mathop{\rm diag}}
\newcommand{\circm}{\mathop{\rm circ}}
\newcommand{\OmegaE}{\Omega_e}

\renewcommand{\epsilon}{\varepsilon}

\renewcommand{\hat}{\widehat}

\newcommand{\bR}{\mathbb{R}}
\newcommand{\bC}{\mathbb{C}}
\newcommand{\bRR}{\mathbb{C}}

\newcommand{\mylenw}{8cm}
\newcommand{\mylenwg}{12cm}

\newcommand{\safeincludegraphics}[2][]{\IfFileExists{#2}{\includegraphics[#1]{#2}}{\fbox{\begin{minipage}[c][45mm][c]{0.90\linewidth}
  \centering
  Missing figure: \texttt{\detokenize{#2}}
\end{minipage}}}}

\everymath{\displaystyle}

\title{
FFT-accelerated computation of \\
the Dirichlet-to-Neumann map for \\
 inhomogeneous exterior Helmholtz problems \\
using the method of fundamental solutions
}

\author{Takemi Shigeta\corref{cor1}\fnref{fn1}}
\ead{shigeta@ac.shoyaku.ac.jp, takemishigeta@gmail.com}
\address{Laboratory of Applied Mathematics, Showa Pharmaceutical University\\
3-2-1 Higashi-Tamagawagakuen, Machida, Tokyo 194-8543, Japan}
\cortext[cor1]{Corresponding author}
\fntext[fn1]{Professor}
\begin{document}

\begin{abstract}
An efficient numerical method is proposed for computing the
Dirichlet-to-Neumann (DtN) map associated with the exterior
Dirichlet problem for the two-dimensional Helmholtz equation
with an inhomogeneous term.
The solution in the subdomain exterior to a circular artificial boundary is approximated by the method of
fundamental solutions (MFS).
When the source and collocation points are equally spaced on
concentric circles, the coefficient matrices arising in the
discretization become circulant, which enables the efficient
evaluation of the discrete DtN map by the fast Fourier
transform (FFT).
By incorporating the transparent boundary condition (TBC) defined by the DtN map
into a finite element formulation, the original exterior problem
posed in an unbounded domain is reduced to an equivalent problem
in a bounded computational domain, which can then be solved by
the finite element method (FEM).
Numerical examples show that the proposed MFS-based approach
provides a highly accurate approximation of the DtN map, that the FFT-based algorithm enables the efficient computation of the resulting discrete DtN map, and that the FEM incorporating the TBC defined by
the DtN map yields accurate solutions for exterior Helmholtz
problems posed in unbounded domains with an inhomogeneous term.
\end{abstract}


\maketitle

\section{Introduction}
Exterior boundary value problems for elliptic partial differential equations arise in many applications, including wave propagation and scattering.
A typical example is the exterior inhomogeneous Helmholtz problem, which appears in time-harmonic acoustic radiation when a localized source is situated near an obstacle.
For instance, waves emitted by a vibrating device or a loudspeaker in the vicinity of
an obstacle can be modeled by the Helmholtz equation with a compactly supported inhomogeneous term in the exterior domain.
Since the radiated field propagates into an unbounded domain, an accurate treatment of the radiation condition is essential.

Among numerical methods for partial differential equations, the finite element method (FEM) is one of the most reliable and widely used methods.
However, the FEM is naturally formulated in bounded computational domains and therefore cannot be applied directly to exterior problems posed in unbounded domains. Therefore, in order to solve exterior Helmholtz problems by the FEM, it is necessary to replace the original problem by an equivalent problem posed in a bounded domain.

In this context, the Dirichlet-to-Neumann (DtN) map plays a central role.
The DtN map associates Dirichlet data prescribed on an artificial
 boundary with the corresponding normal derivative on that boundary.
Once the DtN map is available on the artificial boundary, the original exterior problem can be reformulated as an equivalent boundary value problem in a bounded domain equipped with a transparent boundary condition (TBC) on the artificial boundary~\cite{BTL,NN}.
This makes it possible to solve the original problem in an unbounded domain by using the FEM in a bounded computational domain.

A classical approach to exterior boundary value problems is the Dirichlet--Neumann alternating method (D--N method) \cite{DY,LD,Shigeta2024,Yu1}.
This is an iterative domain decomposition method \cite{Lu}, in which an artificial boundary divides the domain into two subdomains: the inhomogeneous equation is solved in the interior subdomain, while the homogeneous equation is solved in the exterior subdomain.
Its convergence rate depends strongly on the relaxation parameter used to update the boundary values on the artificial boundary.

For the two-dimensional Poisson equation, when both the boundary of the original problem and the artificial boundary are concentric circles, the optimal constant relaxation parameter yielding the fastest convergence is known \cite{Yu2}.
The author further derived an optimal variable relaxation parameter, in which the relaxation parameter is chosen optimally at each iteration \cite{Shigeta2024}.
Although the assumption of circular boundaries simplifies the theoretical analysis, the boundary of the original problem is not necessarily circular in general.
For this reason, it is meaningful to consider artificial boundaries of more general shapes \cite{NN}.
Another possible approach is to use conformal mapping to transform the problem to one with a circular boundary \cite{Liu-c}, so that the optimal variable relaxation parameter can be applied.

However, although conformal mapping is effective for the two-dimensional Poisson equation, its application to the Helmholtz equation becomes difficult because the wavenumber is transformed into a spatially varying coefficient.
Moreover, conformal mapping is not available in three dimensions.
These difficulties motivate the development of a direct numerical method that does not rely on iterative relaxation schemes restricted to circular boundaries or on conformal mapping.

In this paper, we develop such a direct numerical method for solving the exterior Dirichlet problem for the two-dimensional inhomogeneous Helmholtz equation. As noted above, by introducing a circular artificial boundary, we represent the effect of the unbounded subdomain outside the artificial boundary through the DtN map and impose the resulting condition as a TBC on the artificial boundary. This reduces the original exterior problem to a bounded domain problem, which can then be solved by the FEM.

To construct the DtN map numerically, we employ the method of fundamental solutions (MFS)~\cite{BB,Bog,Ei,KO,Katsu2,Mat}, a boundary-type meshfree method for partial differential equations.
In the MFS, the solution is approximated by a finite linear combination of fundamental solutions with their singularities, called source points, placed outside the computational domain. The resulting approximation satisfies the governing equation exactly throughout the domain. Therefore, the expansion coefficients are determined solely by enforcing the boundary conditions at collocation points on the boundary.
Since the MFS does not require mesh generation or volume integration, it is simple to implement, and it is well known to yield highly accurate approximations for problems with smooth geometries and smooth solutions.

The MFS employed in the construction of the DtN map belongs to the broader class of
radial basis function (RBF) methods.
Although RBF methods can provide highly accurate and efficient discretizations~\cite{GCK} and, under suitable node configurations, may allow computations that exploit circulant matrix structures~\cite{Liu-c2,Liu-c}, as described later, the present study couples the MFS-based DtN map with the FEM because the FEM provides a well-established variational framework for inhomogeneous boundary value problems and complex geometries. This coupling also demonstrates that the proposed method for computing the DtN map can be incorporated naturally into a standard computational framework, although it is not intrinsically restricted to the FEM.
A fully meshfree formulation would require additional choices concerning the basis functions, node distribution, shape parameters, and stabilization. A systematic investigation of these issues would constitute a separate study and is beyond the scope of the present paper.

For a circular artificial boundary, the exact DtN map may alternatively be expressed in terms of a Fourier series~\cite{DY}. This analytical construction is highly effective, but it relies on separation of variables and an explicit modal representation of the exterior solution. Boundary integral equations also provide a standard framework for constructing the DtN map~\cite{LD}, although their numerical implementation generally involves singular or nearly singular boundary integrals. 

The MFS is adopted here because, as noted above, it satisfies the exterior Helmholtz equation and the radiation condition without requiring meshing of the exterior domain or singular integration. In addition, it directly constructs a discrete DtN map and permits FFT-based acceleration under rotational symmetry, as described below. Although it may suffer from ill-conditioning depending on the locations of the source points, it remains well suited to the present purpose.

The combination of the MFS with  the fast Fourier transform (FFT) has long been used to solve a wide variety of boundary value problems efficiently in radially symmetric domains. Such methods are commonly referred to as matrix decomposition algorithms (MDAs); see, for example, Chapter~9 of \cite{CCK}.
In the present study, an MDA is employed to compute the discrete DtN map efficiently by exploiting the circulant structure of the resulting matrices \cite{Davis,Liu-c2,Liu-c,SK}.

In a previous paper by the author, a related FFT-based numerical computation of the DtN map for the homogeneous exterior Helmholtz problem was presented in \cite{Shigeta2025}. In the present paper, we extend this line of research to the inhomogeneous Helmholtz equation. The main points are twofold: first,
an MDA for the MFS is applied to the accurate and efficient construction of the DtN map on an artificial boundary;
second, the resulting TBC defined by the DtN map is incorporated into a finite element formulation, which enables  inhomogeneous exterior Helmholtz problems in unbounded domains to be solved in a bounded computational domain.

The remainder of this paper is organized as follows.
Section~\ref{sec:pf} formulates the exterior boundary value problem for the two-dimensional inhomogeneous Helmholtz equation and introduces a circular artificial boundary together with the DtN map.
Section~\ref{sec:mfs} presents the MFS formulation for approximating the solution outside the artificial boundary, together with the resulting discretization of the DtN map.
Section~\ref{sec:fft} develops an FFT-based algorithm for the efficient computation of the discrete DtN map by exploiting the circulant structure of the matrices.
Section~\ref{sec:fem} describes the finite element discretization with a TBC defined by the DtN map.
Section~\ref{sec:ne} presents numerical examples demonstrating the high accuracy of the MFS-based DtN map, the efficiency of the FFT-based computation, the successful incorporation of the TBC defined by the DtN map into the finite element formulation for a homogeneous exterior problem, and the applicability of the resulting method to inhomogeneous exterior problems, including a physically motivated example.
Finally, Section~\ref{sec:conc} concludes the paper.

\section{Problem formulation}\label{sec:pf}
Throughout this paper, we identify $\bR^2$ with $\bC$ via $(x,y)\mapsto x+iy$, and write $z=x+iy$, where $i:=\sqrt{-1}$ denotes the imaginary unit.
Let $R\in C[0,2\pi)$ be a given function.
We consider the simple closed curve
\[
\Gamma:=\{R(\theta)e^{i\theta}:0\leq\theta<2\pi\}
\]
and the simply connected bounded domain $\Omega_0\subset\bRR$
enclosed by $\Gamma$.
Without loss of generality we assume that $0\in\Omega_0$.
We consider the following exterior Dirichlet problem
for the inhomogeneous Helmholtz equation, in which
the unknown function
$u\in H^1_{\rm loc}(\bRR\setminus\overline{\Omega_0})$
is to be determined:
\begin{subequations}\label{eq:orgprob}
\begin{alignat}{7}
-(\Delta+ \kappa^2) u&=f & \quad & \text{in} & \quad &
\bRR\setminus \overline{\Omega_0}, \\
u&=g & & \text{on} & & \Gamma,\\\frac{\partial u}{\partial r}-i\kappa u
&=o(r^{-1/2}) & & \text{as} & & r:=|z|\to\infty ,\label{eq:src}
\end{alignat}
\end{subequations}
where $\Delta:=\partial^2/\partial x^2+\partial^2/\partial y^2$
denotes the Laplace operator,
$\kappa>0$ a given wavenumber,
$f\in L^2(\bRR\setminus\overline{\Omega_0})$
a given inhomogeneous term, and
$g\in H^{1/2}(\Gamma)$ the prescribed Dirichlet data.
Equation \eqref{eq:src} is the Sommerfeld radiation condition, which ensures that the wave propagates outward toward infinity.
We assume that $f\in L^2(\bRR\setminus\overline{\Omega_0})$ has a continuous representative. The support of this representative is then defined by
\[
\supp f:=\overline{\{z\in
\bRR\setminus\overline{\Omega_0}:f(z)\neq 0\}}.
\]
We further assume that $\supp f$ is bounded.
The boundedness of $\supp f$ plays an essential role in the subsequent formulation.
For example, this assumption is satisfied when the wave field is generated by a localized acoustic source, such as a loudspeaker placed near an obstacle.

\begin{figure}[H]
\centering
\includegraphics[width=\mylenw]{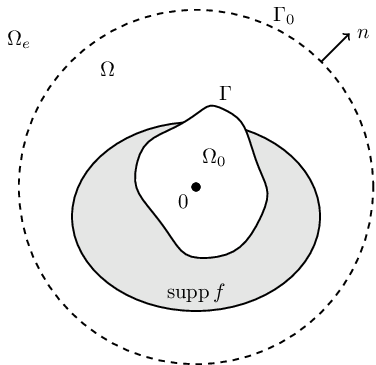}
\caption{Domain decomposition by introducing an artificial boundary}
\label{fig:domain}
\end{figure}

We introduce a circular artificial boundary
$\Gamma_0:=\{R_0e^{i\theta}:0\leq\theta<2\pi\}$
such that $R_0>R(\theta)$ for all $0\leq\theta<2\pi$
and $\Gamma_0$ encloses $\supp f$, and thereby
 decompose the unbounded domain
$\bRR\setminus\overline{\Omega_0}$
into the bounded subdomain
$\Omega$ containing $\supp f$
and the remaining unbounded subdomain
$\OmegaE$.
Note that $f=0$ in $\OmegaE$.
Let $v$ denote the restriction of the unknown function $u$ to $\OmegaE$.
For simplicity, we continue to denote the restriction of $u$ to $\Omega$ by $u$.
Then the following coupling conditions hold on
the artificial boundary $\Gamma_0$:
\begin{equation}
u=v,\quad
\frac{\partial u}{\partial n}
=
\frac{\partial v}{\partial n}
\quad \text{on}\quad \Gamma_0 ,
\label{eq:cnct}
\end{equation}
where $n=n(z):=z/R_0$ denotes the outward unit normal vector on
$\Gamma_0$ with respect to the domain $\Omega$
(see Figure~\ref{fig:domain}).

The Dirichlet-to-Neumann (DtN) map
$\Lambda:H^{1/2}(\Gamma_0)\to H^{-1/2}(\Gamma_0)$
is defined by
\begin{equation}
\Lambda\lambda:=\frac{\partial v}{\partial n},
\qquad \lambda\in H^{1/2}(\Gamma_0),
\label{eq:dnmap}
\end{equation}
where $v\in H^1_{\rm loc}(\OmegaE)$ denotes the solution of
the following exterior Dirichlet problem for the
Helmholtz equation:
\begin{subequations}
\label{eq:extsub}
\begin{alignat}{7}
(\Delta +\kappa^2) v&=0 & \quad & \text{in} & \quad &
\OmegaE, \label{eq:pdeE} \\
v&=\lambda & & \text{on} & & \Gamma_0, \label{eq:bcE} \\
\frac{\partial v}{\partial r}-i\kappa v
&=o(r^{-1/2}) & & \text{as} & & r\to\infty .\label{eq:srcE}
\end{alignat}
\end{subequations}
For the Dirichlet data $\lambda$ in  the Dirichlet condition \eqref{eq:bcE}, we denote the corresponding Neumann data by
$\mu:=\partial v/\partial n$ on $\Gamma_0$.
Here, $\lambda$ and $\mu$ are used as local notation for the Dirichlet and Neumann data on $\Gamma_0$, respectively.
By the definition of the DtN map, we then have
\begin{equation}
\Lambda\lambda=\mu. \label{eq:Ld=n}
\end{equation}In the original problem, however, $\lambda$ is not independently prescribed on $\Gamma_0$, but represents the trace of the solution $u$ from the interior domain:
\[
\lambda=u|_{\Gamma_0}.
\]
Thus, $\lambda$ is an unknown boundary quantity.

Then, from the coupling condition \eqref{eq:cnct}
and the DtN map \eqref{eq:dnmap},
the original problem \eqref{eq:orgprob}
can be reformulated equivalently as the following
mixed boundary value problem posed in the
doubly connected bounded domain,
in which
the unknown function
$u\in H^1(\Omega)$
is to be determined:\begin{subequations}
\label{eq:intsub}
\begin{alignat}{7}
-(\Delta +\kappa^2) u&=f & \quad & \text{in} & \quad & \Omega,\label{eq:fempde}\\
u&=g & & \text{on} & & \Gamma, \label{eq:fembc} \\
\frac{\partial u}{\partial n}
-\Lambda u&=0 & & \text{on} & & \Gamma_0.\label{eq:tbc}
\end{alignat}
\end{subequations}
Equation \eqref{eq:tbc} is called the transparent boundary condition (TBC), which
allows outgoing waves to pass through the artificial boundary without generating spurious reflections, so that the solution inside the bounded domain agrees with the solution of the original exterior problem~\cite{BTL}.

\section{Method of fundamental solutions}\label{sec:mfs}
\subsection{Discretization of the solution}

For $\zeta\in\bC$, the fundamental solution $\Phi_\kappa(\cdot,\zeta)$
to the Helmholtz operator $-(\Delta+\kappa^2)$ is defined by
\[
-(\Delta+\kappa^2)\Phi_\kappa(z,\zeta)=\delta(z-\zeta)
\]
in the sense of distributions, and is given by
\[\Phi_\kappa(z,\zeta):=\frac{i}{4}H_0^{(1)}(\kappa |z-\zeta|),
\qquad z,\zeta\in\mathbb{C},\ z\neq \zeta,
\] where $\delta(z-\zeta)$ denotes the Dirac delta distribution concentrated at
the point $\zeta$,
and $H_\nu^{(1)}$ the Hankel function of the first kind of order $\nu$.

We now briefly show that the fundamental solution $\Phi_\kappa(,\cdot,\zeta)$ satisfies the Sommerfeld radiation condition \eqref{eq:srcE}.
Fix a source point $\zeta$, and write $z=r\widehat z$, where
$r=|z|$ and $|\widehat z|=1$.
As $r\to\infty$, we have
\[
\begin{aligned}
|z-\zeta|
&=
|r\widehat z-\zeta|
=
r-\Re\!\left(\widehat z\,\overline{\zeta}\right)+O(r^{-1})
=
r+O(1).
\end{aligned}
\]
Moreover, differentiating with respect to $r$ while keeping
$\widehat z$ fixed gives
\[
\begin{aligned}
\frac{\partial}{\partial r}|z-\zeta|
&=
1+O(r^{-2}).
\end{aligned}
\]
Set $t:=\kappa|z-\zeta|$.
Since $\kappa>0$ and $|z-\zeta|=r+O(1)$, it follows that $t\to\infty$ as
$r\to\infty$. The large-argument asymptotic expansions
\[
H_0^{(1)}(t)
=
\sqrt{\frac{2}{\pi}}\,
e^{i(t-\pi/4)}
\left(
t^{-1/2}+O(t^{-3/2})
\right)
\]
and
\[
H_1^{(1)}(t)
=
-i\sqrt{\frac{2}{\pi}}\,
e^{i(t-\pi/4)}
\left(
t^{-1/2}+O(t^{-3/2})
\right)
\]
together with
\[
H_0^{(1)\prime}(t)=-H_1^{(1)}(t)
\]
yield
\[
H_0^{(1)\prime}(t)-iH_0^{(1)}(t)
=
O(t^{-3/2}).
\]
The chain rule then gives
\[
\begin{aligned}
\frac{\partial}{\partial r}\Phi_\kappa(z,\zeta)
-i\kappa\Phi_\kappa(z,\zeta)
&=
\frac{i\kappa}{4}
\left[
H_0^{(1)\prime}(t)
\frac{\partial}{\partial r}|z-\zeta|
-iH_0^{(1)}(t)
\right]\\
&=
\frac{i\kappa}{4}
\left[
H_0^{(1)\prime}(t)-iH_0^{(1)}(t)
\right]
+
\frac{i\kappa}{4}
H_0^{(1)\prime}(t)
\left(
\frac{\partial}{\partial r}|z-\zeta|-1
\right) \\
&=O(t^{-3/2})+O(t^{-1/2})O(r^{-2})
=o(r^{-1/2})
\qquad (r\to\infty).
\end{aligned}
\]
Thus, $\Phi_\kappa(\,\cdot\,,\zeta)$ satisfies the Sommerfeld
radiation condition \eqref{eq:srcE}.

\begin{figure}[H]
\centering
\includegraphics[width=\mylenw]{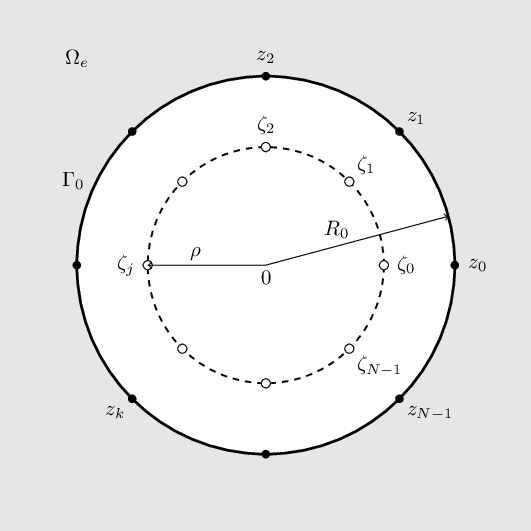}
\caption{Arrangement of source and collocation points}
\label{fig:MFS}
\end{figure}

We choose $N$ distinct source points $\{\zeta_j\}_{j=0}^{N-1}\subset\bRR\setminus\overline{\OmegaE}$ equally spaced on a fictitious circle centered at the origin with radius $\rho$ $(0<\rho<R_0)$, as illustrated in Figure~\ref{fig:MFS}, by setting
\[
\theta_j:=\frac{2\pi j}N,\qquad
\omega:=e^{2\pi i/N},\qquad
\zeta_j:=\rho \omega^j=\rho e^{i\theta_j},
\qquad j=0,1,\ldots,N-1 ,
\]
where $\omega$ is the primitive $N$th root of unity.
Here and throughout the paper, Roman letters are used for the original boundary and its points, while Greek letters are reserved for the fictitious circle and source points.
Then the solution of the exterior Dirichlet problem \eqref{eq:extsub} is discretized by the method of fundamental solutions (MFS)~\cite{BB,Bog,Ei,KO,Katsu2,Mat} as
\begin{equation}
v(z)\approx v^{(N)}(z):=\sum_{j=0}^{N-1}\alpha_j\Phi_\kappa(z,\zeta_j),
\qquad z\in\overline{\OmegaE}:=\OmegaE\cup\Gamma_0,
\label{eq:u_mfs}
\end{equation}
where $\{\alpha_j\}_{j=0}^{N-1}\subset\mathbb{C}$ are unknown expansion coefficients to be determined later.

Since the functions $\{\Phi_\kappa(\cdot,\zeta_j)\}_{j=0}^{N-1}$,
which are used as the MFS basis functions, satisfy the Helmholtz equation \eqref{eq:pdeE} and the Sommerfeld radiation condition \eqref{eq:srcE}, the approximate solution \eqref{eq:u_mfs} also satisfies \eqref{eq:pdeE} and \eqref{eq:srcE}.
Therefore, it remains only to impose the boundary condition \eqref{eq:bcE}.

Next,
we choose $N$ distinct collocation points $\{z_k\}_{k=0}^{N-1}\subset\Gamma_0$ equally spaced along the circular boundary $\Gamma_0$, defined by
\[ z_k:=R_0\omega^k=R_0e^{i\theta_k},
\qquad k=0,1,\ldots,N-1. \]
The arrangement of the collocation points on $\Gamma_0$ is also shown in Figure~\ref{fig:MFS}.
Substituting $z=z_k$ into the Dirichlet condition \eqref{eq:bcE}, we obtain the following linear system for $\{\alpha_j\}_{j=0}^{N-1}\subset\mathbb{C}$:
\[
v^{(N)}(z_k)=\sum_{j=0}^{N-1}\alpha_j
\Phi_\kappa(z_k,\zeta_j)=\lambda(\theta_k),
\qquad k=0,1,\ldots,N-1,
\]
or, in matrix-vector form,
\begin{equation}
C_0\bm\alpha=\bm \lambda ,
\label{eq:C1ag}
\end{equation}
where the matrix $C_0\in\mathbb{C}^{N\times N}$ and
the vectors $\bm \alpha,\bm \lambda\in\mathbb{C}^N$ are defined by
\[
(C_0)_{k,j}:=\Phi_\kappa(z_k,\zeta_j),
\quad
(\bm\alpha)_j:=\alpha_j,
\quad
(\bm\lambda)_k:=\lambda(\theta_k),
\quad
k,j=0,1,\ldots,N-1.
\]

If $C_0$ is nonsingular, then the unknown expansion coefficients are uniquely determined as
\begin{equation}
\bm\alpha=C_0^{-1}\bm \lambda .
\label{eq:aC1g}
\end{equation}
Substituting \eqref{eq:aC1g} into \eqref{eq:u_mfs}, we obtain an approximate solution to \eqref{eq:extsub}
throughout $\overline{\OmegaE}$.

\subsection{Discretization of the DtN map}

Taking the normal derivative of \eqref{eq:u_mfs} and substituting $z=z_k$,
we obtain the approximation of the normal derivative of $v$ at the collocation points on $\Gamma_0$ as follows:
\[
\frac{\partial v^{(N)}}{\partial n}(z_{k})
=\sum_{j=0}^{N-1}\alpha_j
\frac{\partial \Phi_\kappa}{\partial n}(z_k,\zeta_j),
\qquad k=0,1,\ldots,N-1,
\]
or, in matrix-vector form,
\begin{equation}
\bm \mu=C_1\bm\alpha ,
\label{eq:c2}
\end{equation}
where the vector $\bm \mu\in\mathbb{C}^N$ and the matrix
$C_1\in\mathbb{C}^{N\times N}$ are defined by
\[
(\bm \mu)_k:=\frac{\partial v^{(N)}}{\partial n}(z_k),
\qquad
(C_1)_{k,j}:=\frac{\partial \Phi_\kappa}{\partial n}(z_k,\zeta_j),
\qquad k,j=0,1,\ldots,N-1.
\]

The matrices $C_l$ $(l=0,1)$ can be interpreted as matrices whose
entries are given by the $l$th derivatives of the fundamental solution.
The normal derivative of the function $\Phi_\kappa(\cdot,\zeta_j)$ is explicitly expressed as
\[
\frac{\partial \Phi_\kappa}{\partial n}(z,\zeta_j)
=-\frac{i\kappa}{4}
H_1^{(1)}(\kappa|z-\zeta_j|)
\frac{
(z-\zeta_j)\cdot{n(z)}}
{|z-\zeta_j|},
\]
where we define the Euclidean inner product of
$z_1=x_1+iy_1$ and $z_2=x_2+iy_2$ by
\[
z_1\cdot z_2 := \Re\!\left(z_1\overline{z_2}\right) = x_1x_2+y_1y_2 .
\]

Substituting \eqref{eq:aC1g} into \eqref{eq:c2}, we obtain
\begin{equation}
C_1C_0^{-1}\bm \lambda=\bm \mu .
\label{eq:c2c1}
\end{equation}
Equation \eqref{eq:c2c1} corresponds to a discretization of
\eqref{eq:Ld=n}.
Therefore, the matrix
\begin{equation}
\Lambda^{(N)}:=C_1C_0^{-1}\label{eq:ddtn}
\end{equation}
can be regarded as a discrete approximation of the DtN map $\Lambda$.

\section{Efficient computation using the FFT}\label{sec:fft}

In this section we describe an efficient method for computing
the discrete DtN map.
The key observation is that the matrices arising from the
MFS discretization become circulant when both the source
points and collocation points are equally spaced on circles.
This structure enables the discrete DtN map to be computed efficiently by the FFT~\cite{Davis,Liu-c2,Liu-c,SK}.

Although the properties of circulant matrices are standard, we state the relevant details here to fix the conventions used throughout this paper. In particular, the diagonalization formula depends on whether the circulant matrix is defined in terms of its first column or its first row. Moreover, when the Fourier matrix is left unnormalized, the eigenvalues coincide exactly with the components of the discrete Fourier transform of the first column.
When the normalized Fourier matrix is used, the eigenvalues are $\sqrt{N}$ times the components of the corresponding transform.

\subsection{Discrete Fourier Transform}
Let $\bm c=(c_0,c_1,\dots,c_{N-1})^T\in\mathbb{C}^N$.
The discrete Fourier transform (DFT) of $\bm c$ is defined by
\begin{equation}
\hat{c}_k
:=
\sum_{j=0}^{N-1} c_j \omega^{jk},
\qquad k=0,1,\dots,N-1,\label{eq:dft}
\end{equation}
where $\omega= e^{2\pi i/N}$ as defined above.

The DFT \eqref{eq:dft} can be written in the matrix form:
\[
\hat{\bm c} = {\cal F}(\bm c):=W \bm c
\]
with $\hat{\bm c}=(\hat c_0,\hat c_1,\dots,\hat c_{N-1})^T\in\mathbb{C}^N$,
where $W\in\mathbb{C}^{N\times N}$ is the
DFT matrix, also referred to as the unnormalized Fourier matrix, defined by
\[
(W)_{k,j} := \omega^{jk},
\qquad k,j=0,1,\dots,N-1.
\]
Then $W$ satisfies
\[
W^*W=WW^*=NI,
\]
where $W^*:=\overline W^T$ denotes the conjugate transpose of $W$, and $I$  the identity matrix of order $N$.
The inverse DFT (IDFT) is given by
\[
c_j
:=\frac 1N
\sum_{k=0}^{N-1} \hat c_k \omega^{-jk},
\qquad j=0,1,\dots,N-1,
\]
or, in matrix-vector form,
\[
\bm c = {\cal F}^{-1}(\hat{\bm c})
:=W^{-1}\hat{\bm c}=\frac1NW^* \hat{\bm c}.
\]

\begin{remark}
With the above definition of the DFT,
$W$ is the DFT matrix, whereas
$F:=W/\sqrt N$ is the corresponding unitary Fourier matrix,
also referred to as the normalized Fourier matrix, satisfying $F^*F=FF^*=I$.
\end{remark}

\subsection{Circulant matrices}
A matrix $C\in\mathbb{C}^{N\times N}$
is called a circulant matrix if each column vector is obtained from the previous one by a cyclic downward shift.
Equivalently, if the first column of $C$ is $\bm c=(c_0,c_1,\dots,c_{N-1})^T$, then
the $(k,j)$ entry of $C$ is given by
\[
(C)_{k,j}=c_{(k-j) \bmod N},\qquad k,j=0,1,\ldots,N-1,
\]
where the indices are understood modulo $N$. In particular,
$c_{-j}=c_{N-j}$ for $j=1,2,\dots,N-1$.
In this case, we write
\[ C=\circm(\bm c),\qquad
\circm(\bm c):=\begin{pmatrix}
c_0 & c_{N-1} & c_{N-2} & \cdots &c_1 \\
c_1 & c_{0} & c_{N-1} & \cdots &c_2 \\
c_2 & c_{1} & c_{0} & \cdots &c_3 \\
\vdots & \vdots & \vdots & \ddots &\vdots\\
c_{N-1} & c_{N-2} & c_{N-3} & \cdots &c_0
\end{pmatrix}.
 \]

We consider the eigenvalue problem for the circulant matrix $C$:
\[ C\bm w_k=\sigma_k\bm w_k,\qquad \bm w_k\neq\bm 0,
\qquad k=0,1,\ldots,N-1. \]
It is well known that the eigenvalues of $C$ are given by
\begin{equation}
 \sigma_k = \sum_{j=0}^{N-1} c_j \omega^{jk},\qquad
k=0,1,\ldots,N-1, \label{eq:chat}
\end{equation}
which coincides exactly with the DFT of $\bm c$, given in \eqref{eq:dft}, and hence $\sigma_k=\hat{c}_k$.
For each $k=0,1,\ldots,N-1$, an eigenvector $\bm w_k\in\bC^N$ of $C$ corresponding to
the eigenvalue  $\sigma_k$
is given by
\[
(\bm w_k)_j:=\omega^{-jk},
\qquad
j=0,1,\ldots,N-1. \]
Indeed, this can be verified directly as follows.
For each $l=0,1,\ldots,N-1$,
\begin{align*}
(C\bm w_k)_l
&=\sum_{j=0}^{N-1} (C)_{l,j}(\bm w_k)_j
=
\sum_{j=0}^{N-1} c_{l-j}\omega^{-jk}
=\sum_{m=0}^{N-1} c_m \omega^{-(l-m)k}\\
&=
\omega^{-lk}\sum_{m=0}^{N-1} c_m \omega^{mk}
=\sigma_k\,\omega^{-lk}
=
\sigma_k (\bm w_k)_l ,
\end{align*}
where the indices are understood modulo $N$.
Hence, defining the matrix
\begin{equation}
U:=\begin{pmatrix}
\bm w_0 & \bm w_1 & \cdots \bm w_{N-1}\end{pmatrix},
\label{eq:def_u}
\end{equation}
we obtain the diagonalization
\[
U^{-1}CU=\diag(\hat{\bm c}),
\]
or equivalently,
\[
C=U\diag(\hat{\bm c})U^{-1},
\]
where $\diag(\widehat{\bm c})$ denotes the diagonal matrix whose
diagonal entries are the components of $\widehat{\bm c}$.
Since $U=W^*$, where $W$ is the DFT matrix, this can be rewritten as
\begin{equation}
C=\frac1N W^* \diag(\hat{\bm c}) W.
\label{eq:diagC}
\end{equation}
\begin{remark}
If the DFT and the IDFT are defined by means of the Fourier matrix $F=W/\sqrt N$ as follows:
\[
\hat{\bm c}=\mathcal{F}(\bm c):=F\bm c,\qquad
\bm c=\mathcal{F}^{-1}(\hat{\bm c}):=F^{-1}\hat{\bm c}=F^*\hat{\bm c},
\]
then the eigenvalues satisfy $\sigma_k=\sqrt{N}\,\hat c_k$, and hence
\[
C=F^*\diag(\sqrt{N}\,\hat{\bm c})F.
\]
\end{remark}
\begin{remark}\label{rem:Crow}
Let a circulant matrix $C$ be defined by
\[
(C)_{k,j}=c_{(j-k)\bmod N},\qquad k,j=0,1,\ldots,N-1,
\]
where $\bm c^T=(c_0,c_1,\dots,c_{N-1})$ is the first row of $C$.
Equivalently, $C$ can be written as
\[
C=
\begin{pmatrix}
c_0 & c_1 & c_2 & \cdots & c_{N-1}\\
c_{N-1} & c_0 & c_1 & \cdots & c_{N-2}\\
c_{N-2} & c_{N-1} & c_0 & \cdots & c_{N-3}\\
\vdots & \vdots & \vdots & \ddots & \vdots\\
c_1 & c_2 & c_3 & \cdots & c_0
\end{pmatrix}.
\]
Then, for each $k=0,1,\ldots,N-1$, an eigenvector $\bm w_k\in\bC^N$ of $C$ corresponding to
the eigenvalue  $\sigma_k$ in \eqref{eq:chat}
is given by
\[
(\bm w_k)_j:=\omega^{jk},
\qquad
j=0,1,\ldots,N-1. \]
If we define $U$ in the same way as in \eqref{eq:def_u}, then $U=W$.
Therefore,
the circulant matrix $C$ is diagonalized as
\[
C=\frac1N W \diag(\hat{\bm c}) W^*.
\]
Thus, the expression for the diagonalization of a circulant matrix depends on whether the matrix is viewed with respect to its first column or its first row.
\end{remark}
\subsection{Application to the discrete DtN map}
As described in the previous subsection, the matrices $C_0$
and $C_1$ appearing in the MFS discretization are circulant
because both the source points and collocation points are
equally spaced on circles.
Let $\bm c_l$ $(l=0,1)$ denote the first column of $C_l$.
Applying the diagonalization formula for circulant matrices,
we obtain
\[
C_l = \frac1NW^* \operatorname{diag}(\hat{\bm c}_l) W,
\qquad l=0,1,
\]
where
\[
\hat{\bm c}_l ={\cal F}(\bm c_l)= W \bm c_l, \qquad
\bm c_l ={\cal F}^{-1}(\hat{\bm c}_l)= \frac1NW^* \hat{\bm c}_l,
\qquad l=0,1.
\]
Hence, the discrete DtN map applied to $\bm\lambda$ can be written as
\begin{align}
\bm\mu=
\Lambda^{(N)}\bm\lambda
&= C_1 C_0^{-1}\bm\lambda
=\frac1NW^*
\diag(\hat{\bm c}_1)\diag(\hat{\bm c}_0)^{-1}
W\bm\lambda \notag \\
&={\cal F}^{-1}(({\cal F}({\bm c}_1)
\oslash {\cal F}({\bm c}_0)
\odot {\cal F}({\bm \lambda})),
\label{eq:c1c0}
\end{align}
where for two vectors $\bm a=(a_0,a_1,\ldots,a_{N-1})^T$ and
$\bm b=(b_0,b_1,\ldots,b_{N-1})^T$ in $\mathbb{C}^N$,
the Hadamard product $\bm a\odot\bm b$ and the Hadamard division
$\bm a\oslash\bm b$ are defined, respectively, by
\begin{align*}
\bm a\odot\bm b &:= (a_0 b_0, a_1 b_1, \ldots, a_{N-1} b_{N-1})^T,\\
\bm a\oslash\bm b &:= \left(\frac{a_0}{b_0},\frac{a_1}{b_1},\ldots,\frac{a_{N-1}}{b_{N-1}}\right)^T,
\qquad b_j\neq 0,\quad j=0,1,\ldots,N-1 .
\end{align*}
\begin{remark}
As noted in Remark~\ref{rem:Crow}, the diagonalization formula for a circulant matrix depends on whether the matrix is defined by its first column or by its first row.
In the present setting, however, the matrices $C_l$ $(l=0,1)$ are symmetric, so this distinction causes no essential difference.
\end{remark}
\subsection{FFT-based algorithm}
As shown in \eqref{eq:c1c0}, the discrete DtN map can be computed
using the DFT and the IDFT, and therefore can be evaluated efficiently by the fast Fourier transform (FFT) and the inverse FFT (IFFT).
The resulting FFT-based procedure can be summarized as follows.
\begin{description}
\item[Step 1.] Compute the DFTs by the FFT:
\[
\hat{\bm c}_0 := {\cal F}(\bm c_0),\qquad
\hat{\bm  c}_1 := {\cal F}(\bm c_1),\qquad
\hat{\bm \lambda} := {\cal F}(\bm \lambda).
\]
\item[Step 2.] Compute
\[
\hat{\bm r} := \hat{\bm c}_1 \oslash \hat{\bm c}_0 ,
\]
which coincides with the DFT of the first column of
$\Lambda^{(N)}$.
Therefore, the discrete DtN map can be efficiently computed as
\[
\Lambda^{(N)}=\circm\!\left({\cal F}^{-1}(\hat{\bm r})\right).
\]
\item[Step 3.] Compute
\[
\hat{\bm \mu}:=\hat{\bm r} \odot \hat{\bm \lambda}.
\]
\item[Step 4.] Compute the IDFT by the IFFT:
\[
\bm \mu:={\cal F}^{-1}(\hat{\bm \mu}).
\]
\end{description}
It should be noted that the discrete DtN matrix $\Lambda^{(N)}$ is obtained in Step~2. Steps~3 and 4 describe its application to a given vector of Dirichlet data in order to compute the corresponding discrete Neumann data. In the finite element formulation presented in Section~\ref{sec:fem}, only the matrix $\Lambda^{(N)}$ obtained in Step~2 is required.
A straightforward computation of $C_1C_0^{-1}\bm\lambda\,(=\bm\mu)$ requires first solving the dense linear system \eqref{eq:C1ag}, corresponding to \eqref{eq:aC1g}, at a cost of $O(N^3)$ operations, and then evaluating \eqref{eq:c2}, which requires a further $O(N^2)$ operations.
It also requires $O(N^2)$ memory to store the dense matrices.
Thus, the direct approach involves dense-matrix computations
with high computational and memory costs.
By contrast, the present method exploits the circulant structure of $C_0$ and $C_1$.
If the discrete DtN map itself needs to be represented explicitly as a matrix, both the computational cost and the memory requirement are $O(N^2)$. However, if the objective is only to compute the Neumann data corresponding to given Dirichlet data,
the FFT enables the computation to be carried out in $O(N\log N)$ operations without explicitly solving a dense linear system.
Moreover, since each of the circulant matrices $C_0$ and $C_1$ is completely determined by its first column, it suffices to store only these columns, resulting in an $O(N)$ memory requirement.
Hence, the proposed approach significantly reduces both the computational and memory costs.
In fact,
since \eqref{eq:aC1g} can be written as
\[
\bm \alpha=C_0^{-1}\bm\lambda
=\frac1NW^*\diag(\hat{\bm c}_0)^{-1}W\bm\lambda
={\cal F}^{-1}\bigl({\cal F}(\bm\lambda)\oslash{\cal F}(\bm c_0)\bigr),
\]
the unknown expansion coefficients can be computed efficiently by the FFT, without directly solving the linear system \eqref{eq:C1ag}.
\section{Finite element formulation with the transparent boundary condition}\label{sec:fem}
In this section, we briefly
explain how the transparent boundary condition (TBC) is incorporated into the finite element method (FEM) for the interior problem posed in $\Omega$.
Multiplying the Helmholtz equation \eqref{eq:fempde} by the complex conjugate $\overline{w}$ of a test function $w$
and applying Green's formula, we obtain
\[
\int_{\Omega} \nabla u \cdot \nabla \overline w \, d\Omega
-
\kappa^2 \int_{\Omega} u\overline w \, d\Omega
-
\int_{\Gamma_0} (\Lambda u)\,\overline w\,d\Gamma
=
\int_{\Omega} f \overline w\,d\Omega
\]
for all test functions $w\in H^1(\Omega)$  satisfying $w=0$ on $\Gamma$.
Let $V_h \subset H^1(\Omega)$ be the finite element space consisting of piecewise linear functions.
Let $f_h$ and $g_h$ denote the discrete approximations of $f$ and $g$, respectively.
We seek $u_h \in V_h$ such that
\[
u_h=g_h
\qquad \text{on} \qquad \Gamma,
\]
and
\[
\int_{\Omega} \nabla u_h \cdot \nabla \overline w_h \, d\Omega
-
\kappa^2 \int_{\Omega} u_h \overline w_h \, d\Omega
-
\int_{\Gamma_0} (\Lambda u_h)\,\overline w_h\,d\Gamma
=
\int_{\Omega} f_h \overline w_h\,d\Omega
\]
for all $w_h \in V_h$ satisfying
$w_h=0$ on $\Gamma$.
Let $\{\phi_j\}_{j=1}^M$ be the nodal basis functions of $V_h$, with $M$ denoting the number of finite element nodes.
Then the finite element solution $u_h \in V_h$ is expressed as
\[
u_h
=
\sum_{j=1}^M u_j \phi_j,
\]
where $\bm u=(u_1,\dots,u_M)^\top$ denotes the coefficient vector.
We order the nodal degrees of freedom according to the nodes in $\Omega$, on $\Gamma$, and on $\Gamma_0$.
Accordingly, the coefficient vector is decomposed as
\[
\bm u=
\begin{pmatrix}
\bm u_{\Omega}\\
\bm u_{\Gamma}\\
\bm u_{\Gamma_0}
\end{pmatrix}.
\]
With this ordering, we introduce the block matrix
\[
A=
\begin{pmatrix}
A_{\Omega\Omega} & A_{\Omega\Gamma} & A_{\Omega\Gamma_0}\\
A_{\Gamma\Omega} & A_{\Gamma\Gamma} & A_{\Gamma\Gamma_0}\\
A_{\Gamma_0\Omega} & A_{\Gamma_0\Gamma} & A_{\Gamma_0\Gamma_0}
\end{pmatrix},
\]
whose blocks are defined by
\[
(A_{XY})_{k,j}
=
\int_{\Omega}
\nabla \phi_j^{(Y)} \cdot \nabla \overline{\phi_k^{(X)}}
\,d\Omega
-
\kappa^2
\int_{\Omega}
\phi_j^{(Y)} \overline{\phi_k^{(X)}}
\,d\Omega,
\qquad
X,Y \in \{\Omega,\Gamma,\Gamma_0\}
\]
with $\phi_j^{(Y)}$ and $\phi_k^{(X)}$ denoting the basis functions associated with the nodes
belonging to $Y$ and $X$, respectively.
Similarly, the load vector is decomposed as
\[
\bm f=
\begin{pmatrix}
\bm f_{\Omega}\\
\bm f_{\Gamma}\\
\bm f_{\Gamma_0}
\end{pmatrix}
\]
with
\[
(\bm f_X)_k
=
\int_{\Omega} f_h\,\overline{\phi_k^{(X)}}\,d\Omega,
\qquad
X \in \{\Omega,\Gamma,\Gamma_0\}.
\]
Since the nodes on $\Gamma_0$ coincide with the $N$ equally spaced
collocation points used to construct the discrete DtN map,
the subvector $\bm u_{\Gamma_0}$ consists precisely of the values of
$u_h$ at these points.
For $j=1,\ldots,N$, define
\[
\bm\phi_j^{(\Gamma_0)}
:=
\bigl(
\phi_j^{(\Gamma_0)}(\bm x_1),
\ldots,
\phi_j^{(\Gamma_0)}(\bm x_N)
\bigr)^T, \]
where $\bm x_1,\ldots,\bm x_N$ are the nodes on $\Gamma_0$,
corresponding to the equally spaced points $z_0,\ldots,z_{N-1}$ introduced in Section~\ref{sec:mfs}.
Since $\{\phi_j^{(\Gamma_0)}\}_{j=1}^N$ is the nodal basis associated
with these nodes, we have
\[
\bm\phi_j^{(\Gamma_0)}=\bm e_j,
\]
where $\bm e_j$ denotes the $j$th standard basis vector of $\mathbb C^N$.
Applying the $N$-point trapezoidal rule to the boundary integral involving the DtN map, we obtain
\[
\begin{aligned}
(B_{\mathrm{DtN}})_{k,j}
&:=
\int_{\Gamma_0}
(\Lambda\phi_j^{(\Gamma_0)})
\overline{\phi_k^{(\Gamma_0)}}\,d\Gamma
\\
&\approx
\frac{2\pi R_0}{N}
\sum_{l=1}^{N}
\bigl(
\Lambda^{(N)}
\bm\phi_j^{(\Gamma_0)}
\bigr)_l
\overline{
\bigl(
\bm\phi_k^{(\Gamma_0)}
\bigr)_l
}
\\
&=
\frac{2\pi R_0}{N}
{\bm\phi_k^{(\Gamma_0)}}^{*}
\Lambda^{(N)}
\bm\phi_j^{(\Gamma_0)}
\\
&=
\frac{2\pi R_0}{N}
\bm e_k^{*}
\Lambda^{(N)}
\bm e_j
=
\frac{2\pi R_0}{N}
\bigl(\Lambda^{(N)}\bigr)_{k,j}.
\end{aligned}
\]
Accordingly, the matrix representing the boundary integral involving the DtN map is defined by
\[
B_{\mathrm{DtN}}
:=
\frac{2\pi R_0}{N}\Lambda^{(N)}.
\]
Since the DtN map acts on the trace of the finite element function on $\Gamma_0$, the matrix arising from the boundary integral involving the DtN map contributes only to the $(\Gamma_0,\Gamma_0)$ block of the finite element matrix.
Then the resulting finite element linear system is written as
\[
\begin{pmatrix}
A_{\Omega\Omega} & A_{\Omega\Gamma} & A_{\Omega\Gamma_0}\\
A_{\Gamma\Omega} & A_{\Gamma\Gamma} & A_{\Gamma\Gamma_0}\\
A_{\Gamma_0\Omega} & A_{\Gamma_0\Gamma} & A_{\Gamma_0\Gamma_0}-B_{\mathrm{DtN}}
\end{pmatrix}
\begin{pmatrix}
\bm u_{\Omega}\\
\bm u_{\Gamma}\\
\bm u_{\Gamma_0}
\end{pmatrix}
=
\begin{pmatrix}
\bm f_{\Omega}\\
\bm f_{\Gamma}\\
\bm f_{\Gamma_0}
\end{pmatrix}.
\]
Since the Dirichlet boundary condition is prescribed on $\Gamma$, we have
\[
\bm u_{\Gamma}=\bm g,
\]
where $\bm g$ denotes the vector of nodal values of the prescribed boundary data on $\Gamma$.
Therefore, the unknowns are reduced to $\bm u_{\Omega}$ and $\bm u_{\Gamma_0}$, and the resulting reduced system becomes
\[
\begin{pmatrix}
A_{\Omega\Omega} & A_{\Omega\Gamma_0}\\
A_{\Gamma_0\Omega} & A_{\Gamma_0\Gamma_0}-
B_{\mathrm{DtN}}
\end{pmatrix}
\begin{pmatrix}
\bm u_{\Omega}\\
\bm u_{\Gamma_0}
\end{pmatrix}
=
\begin{pmatrix}
\bm f_{\Omega}\\
\bm f_{\Gamma_0}
\end{pmatrix}
-
\begin{pmatrix}
A_{\Omega\Gamma}\\
A_{\Gamma_0\Gamma}
\end{pmatrix}
\bm g .
\]
This is the final linear system to be solved.
\section{Numerical examples}\label{sec:ne}
\subsection{Accuracy of the discrete DtN map}\label{sec:dtntest}
\begin{figure}[H]\centering
\includegraphics[width=\mylenw]{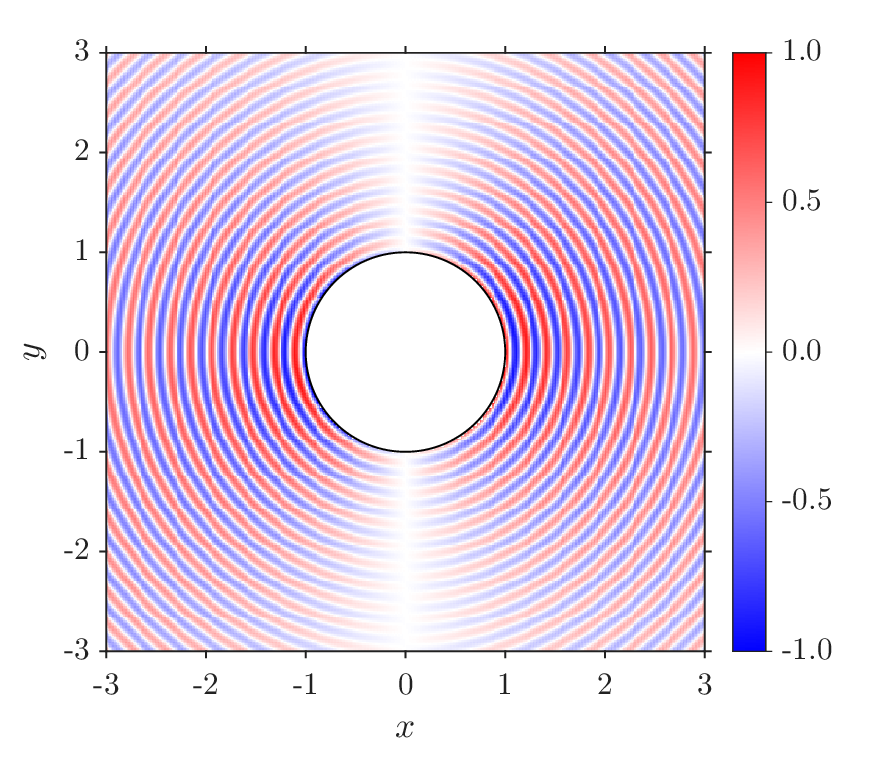}
\caption{Real part of the exact solution, $\Re v$, in $D$}
\label{fig:domexact}
\end{figure}

\begin{figure}[H]\centering
\includegraphics[width=\mylenw]{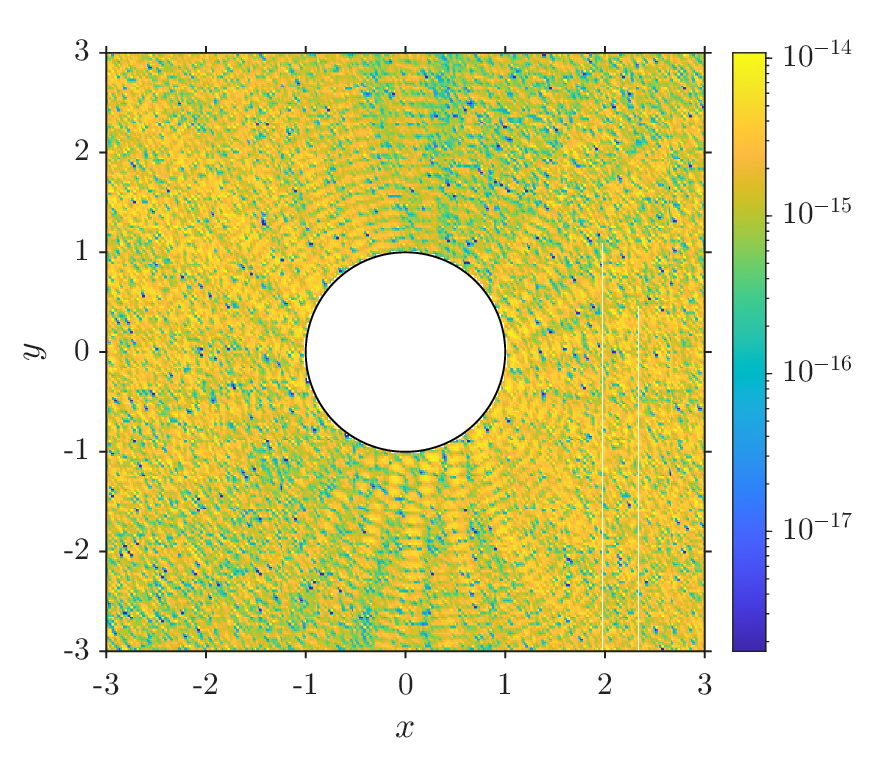}
\caption{Absolute error $\left|\Re v^{(N)}-\Re v\right|$ in $D$}
\label{fig:domerr}
\end{figure}

\begin{figure}[H]\centering
\includegraphics[width=\mylenwg]{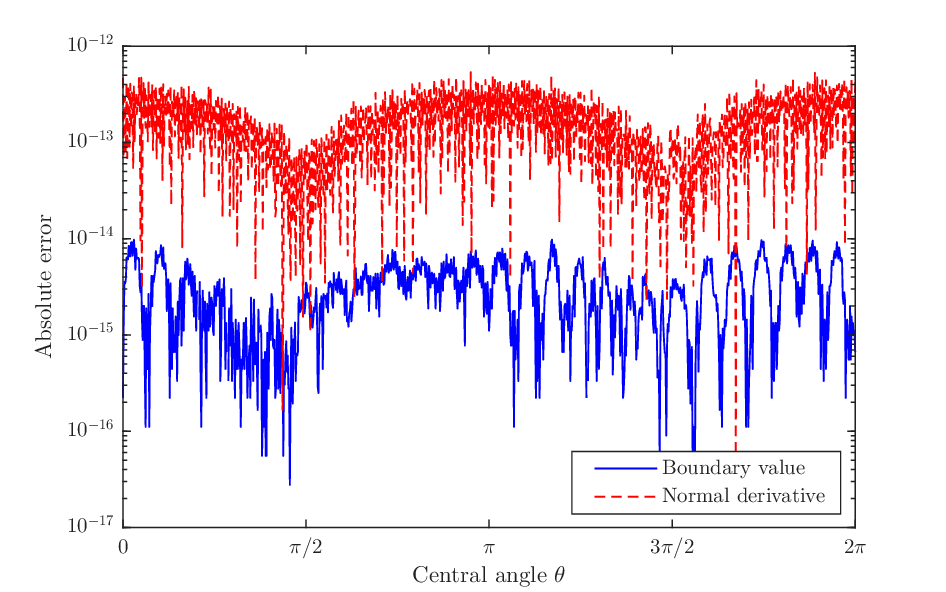}
\caption{Absolute errors in the real part of $v^{(N)}$ and its normal derivative on $\Gamma_0$}
\label{fig:bdyerr}
\end{figure}

For prescribed Dirichlet data $\lambda$, we solve Problem~\eqref{eq:extsub}. We first examine the accuracy of the MFS approximation $v^{(N)}$ to the exact solution $v$ in the exterior subdomain $\OmegaE$ and then assess the accuracy of the corresponding Neumann data $\mu$ on the artificial boundary $\Gamma_0$, thereby verifying the accuracy of the MFS-based discrete DtN map~\cite{Shigeta2025}.
The radius of the circular artificial boundary $\Gamma_0$ is taken as $R_0=1$,
and the wavenumber is set to $\kappa=30$.
The exact solution is given by
\[
v(re^{i\theta})
:=\frac{H_1^{(1)}(\kappa r)}{H_1^{(1)}(\kappa)}\cos\theta ,
\qquad r>1,\quad 0\leq\theta<2\pi.
\]
Figure~\ref{fig:domexact} shows the real part of the exact solution
 in the subdomain
\[
D:=\{z=x+iy\in\bC:-3<x<3,\ -3<y<3,\ |z|>1\},
\]
that is, the domain inside the open square $(-3,3)^2$ and outside the unit circle.
In this case, the Neumann data $\mu$, namely, the normal derivative of $v$ on $\Gamma_0$ is obtained as
\[
\mu(\theta)=\frac{\partial v}{\partial n}(e^{i\theta})
=-\frac{\kappa}{2H_1^{(1)}(\kappa)}
\Bigl(H_0^{(1)}(\kappa)-H_2^{(1)}(\kappa)\Bigr)\cos\theta ,
\qquad 0\leq\theta<2\pi.
\]
Taking the exact solution $v$ as the unknown function,
we solve Problem~\eqref{eq:extsub} with the Dirichlet data
$\lambda(\theta):=\cos\theta$.
\begin{figure}[H]\centering
\includegraphics[width=\mylenwg]
{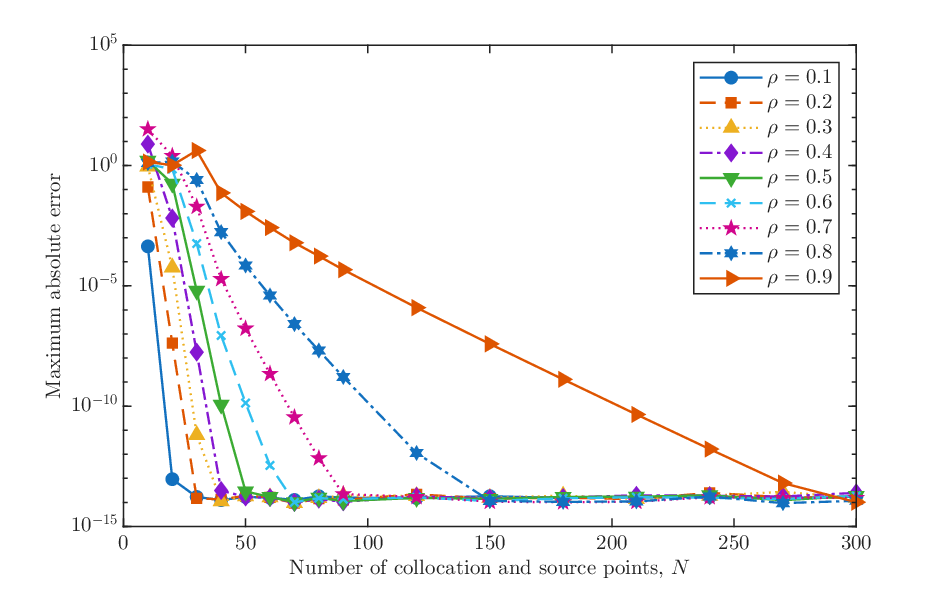}
\caption{Maximum absolute error in $\overline D$ versus $N$ for various values of $\rho$}
\label{fig:err_N}
\end{figure}
Figure~\ref{fig:err_N} shows a semi-logarithmic plot in which the
horizontal axis represents the number of collocation and source points $N$
and the vertical axis represents the maximum absolute error
over $\overline D$:
\begin{align*}
 \|\Re v^{(N)} - \Re v\|_{C(\overline D)}
&:=\max_{z\in\overline D}
\left|\Re v^{(N)}(z)-\Re v(z)\right| \\
&\approx
\max_{1\leq j\leq 9324}
\left|\Re v^{(N)}(z_j)-\Re v(z_j)\right|.
\end{align*}
Here, the points $\{z_j\}_{j=1}^{9324}$ denote the error-evaluation points obtained from a uniform $101\times101$ Cartesian grid in $[-3,3]^2$, excluding the points lying inside the unit circle. They are distinct from the boundary collocation points $\{z_k\}_{k=0}^{N-1}$ introduced in Section~\ref{sec:mfs}.
The number of collocation points was chosen in the range
$10 \leq N \leq 90$ with increments of $10$,
and $90 \leq N \leq 300$ with increments of $30$.
Each curve corresponds to a different value of $\rho$,
and the results for $\rho=0.1$ to $\rho=0.9$ are shown in the legend.
The error decreases overall as $N$ increases.
When $\rho$ is small, the error decays rapidly as $N$ increases and reaches machine precision.
As $\rho$ becomes larger, however, the convergence becomes significantly slower.
Nevertheless, before reaching machine precision, the error still shows an approximately exponential decay.
These results indicate that, also for the Helmholtz equation, the error decays exponentially with respect to $N$, as observed for the MFS when applied to the Laplace equation \cite{Katsu2}.
Next, for the parameters $N=300$ and $\rho=0.9$,
Figure~\ref{fig:domerr} shows the absolute error
in the real part of the numerical solution:
\[
\left|\Re v^{(N)}(z)-\Re v(z)\right|, \qquad z\in  D.
\]
It can be observed that the absolute error is less than
approximately $10^{-14}$,
indicating that a highly accurate numerical solution is obtained.
Furthermore, on $\Gamma_0$ ($r=1$, $0\leq\theta<2\pi$),
the absolute error in the real part of the numerical solution
\[
\left|\Re v^{(N)}(e^{i\theta})-\Re v(e^{i\theta})\right|,
\qquad 0\leq\theta<2\pi
\]
and the absolute error in the real part of the normal derivative
\[
\left|
\Re \frac{\partial v^{(N)}(e^{i\theta})}{\partial n}
-
\Re \frac{\partial v(e^{i\theta})}{\partial n}
\right|,
\qquad 0\leq\theta<2\pi
\]
are shown in
Figure~\ref{fig:bdyerr}, respectively.
This result indicates that the discrete DtN map functions properly
and that the normal derivative on the boundary is computed
with sufficient accuracy even for a large wavenumber.
\subsection{Computational efficiency of the discrete DtN map}
\begin{figure}[H]
\centering
\includegraphics[width=\mylenwg]{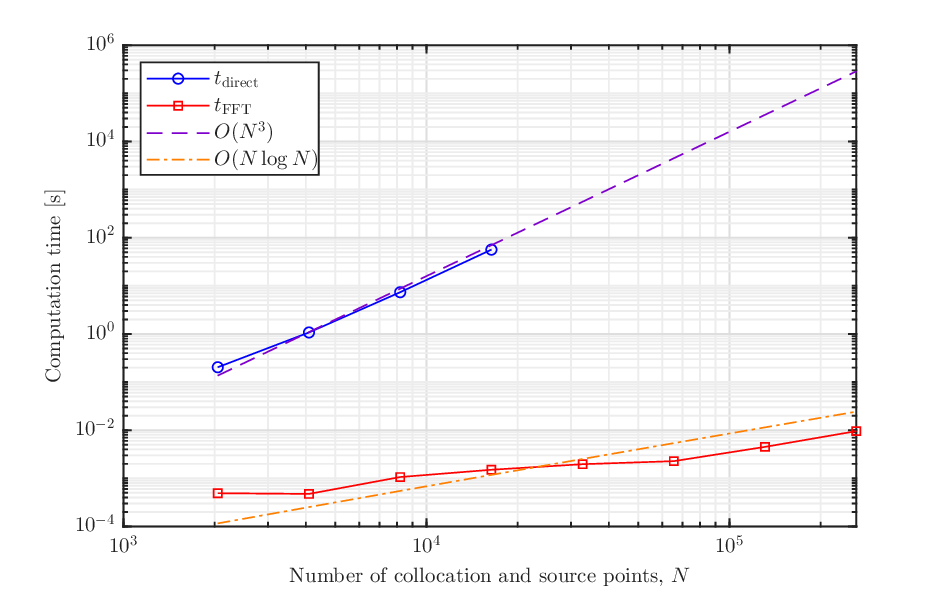}
\caption{Computation time for obtaining the Neumann data from the given Dirichlet data versus the number of collocation and source points $N$}
\label{fig:dtn_time}
\end{figure}

\begin{figure}[H]
\centering
\includegraphics[width=\mylenwg]{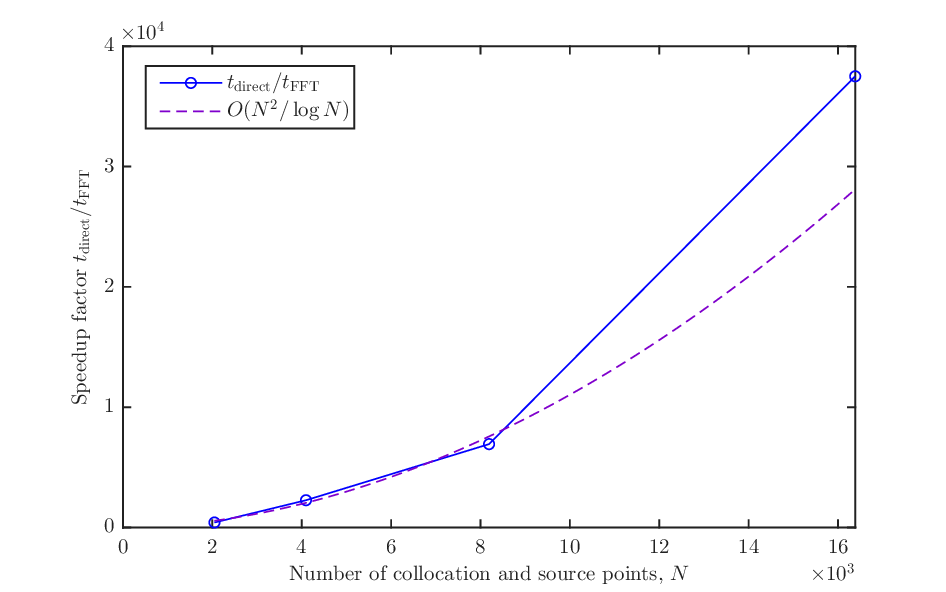}
\caption{Speedup factor of the proposed method relative to the direct method}
\label{fig:dtn_speedup}
\end{figure}

We next examine the computational efficiency of the discrete DtN map.
We use the same setting as in Section~\ref{sec:dtntest}, except that $\kappa=1$ and $\rho=0.99$.
As observed in Figure~\ref{fig:err_N}, choosing $\rho$ closer to the radius of the artificial boundary helps maintain stable numerical accuracy for large values of $N$.
All computation times reported in this subsection were measured in MATLAB R2025b on the same computing environment for both the proposed FFT-based method and the conventional direct method.
The experiments were performed on a PC equipped with
an AMD Ryzen 7 PRO 5750G processor and 16 GB of RAM.
For each value of $N$, the reported computation time is the average of three runs.
Figure~\ref{fig:dtn_time} shows the computation time required to obtain the Neumann data from the given Dirichlet data using the discrete DtN map as a function of $N$.
The results clearly show that the proposed method is significantly faster than the direct method.
As $N$ increases, the computation time $t_{\rm direct}$ for the direct method grows rapidly, whereas the computation time $t_{\rm FFT}$ for the proposed method increases much more mildly.
Moreover, the observed growth is consistent with the expected complexities $O(N^3)$ for the direct method and $O(N\log N)$ for the proposed
method.
This behavior reflects the fact that the direct method requires solving a dense linear system, while the proposed method exploits the circulant structure of the coefficient matrices and evaluates the discrete DtN map efficiently by the FFT.
It is also worth noting that, even for relatively large values of $N$, the proposed method remains faster than the direct method for much smaller values of $N$.
This clearly demonstrates the practical effectiveness of the FFT-based acceleration.
To further illustrate the gain in efficiency, Figure~\ref{fig:dtn_speedup} shows the speedup factor $t_{\rm direct}/t_{\rm FFT}$, defined as the ratio of the computation time of the direct method to that of the proposed method.
The speedup factor increases as $N$ becomes larger, showing that the advantage of the proposed method becomes more pronounced for finer discretizations.
The reference line $O(N^2/\log N)$, derived from the ratio $O(N^3)/O(N\log N)$, is included for comparison.
The measured speedup factor at the largest value of $N$ deviates from the asymptotic trend, possibly because of machine-dependent effects, whereas those at the other values of $N$ closely follow the reference line.
In particular, for large values of $N$, the proposed method achieves a substantial reduction in computational time compared with the direct method.
These results confirm the computational efficiency of the proposed method.
\subsection{Numerical verification of the TBC}
This example demonstrates that the original exterior problem posed in an unbounded domain can be reduced to an equivalent problem in a bounded domain by incorporating the TBC defined by the DtN map into the finite element formulation.
We consider the homogeneous case $f=0$ with a circular boundary $\Gamma$ of radius $R(\theta):=1$.
As an exact solution, we use the single outgoing Fourier--Hankel mode:
\[
u(re^{i\theta})
=
\frac{H_m^{(1)}(\kappa r)}{H_m^{(1)}(\kappa R_0)}e^{im\theta},\qquad r>1,\qquad 0\leq\theta<2\pi,
\]
which satisfies the Helmholtz equation in the exterior domain and the Sommerfeld radiation condition.
The Dirichlet boundary data on $\Gamma$ are taken from this exact solution.
We choose the mode number $m=3$.
\begin{figure}[H]\centering
\includegraphics[width=\mylenw]{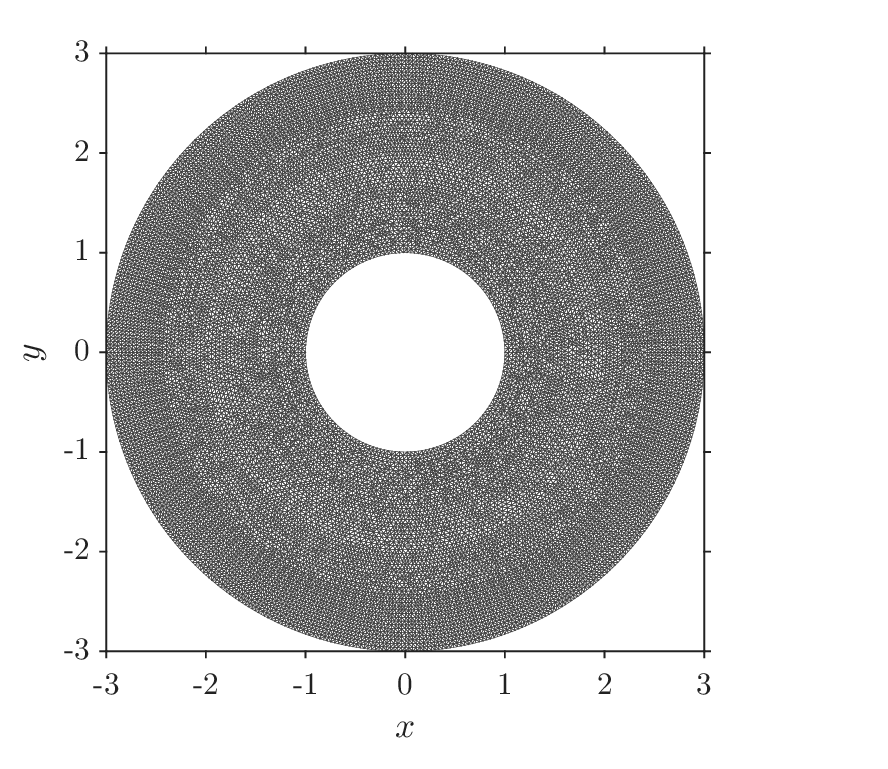}
\caption{Finite element mesh
(19,320 nodes and 37,974 elements)}
\label{fig:tbc_mesh}
\end{figure}
\begin{figure}[H]\centering
\includegraphics[width=\mylenw]{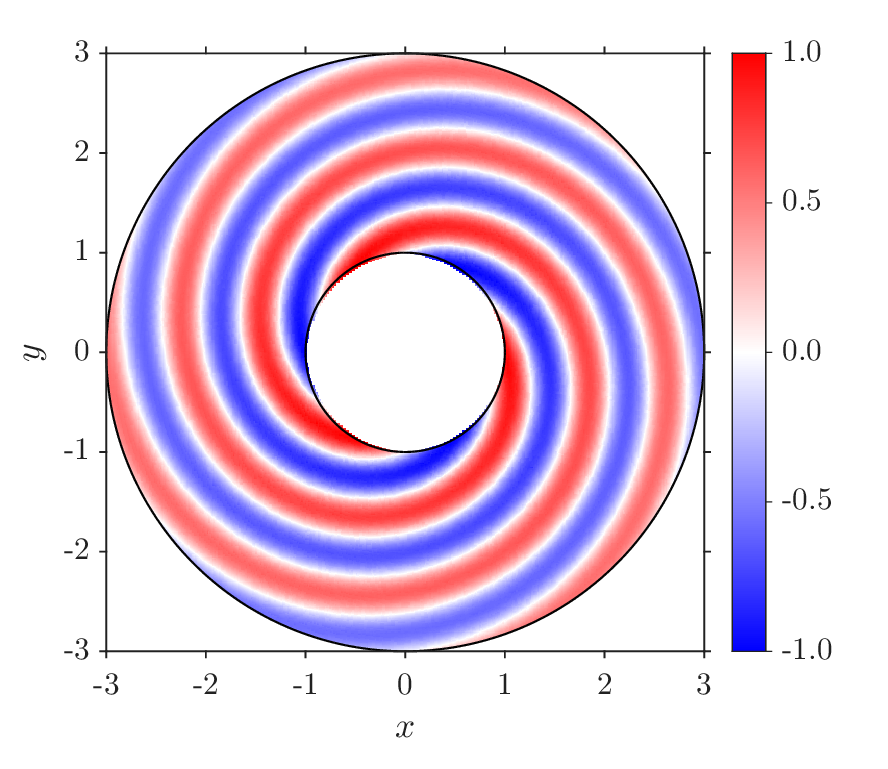}
\caption{Real part of the exact solution, $\Re u$}
\label{fig:tbc_uexRe}
\end{figure}
\begin{figure}[H]\centering
\includegraphics[width=\mylenw]{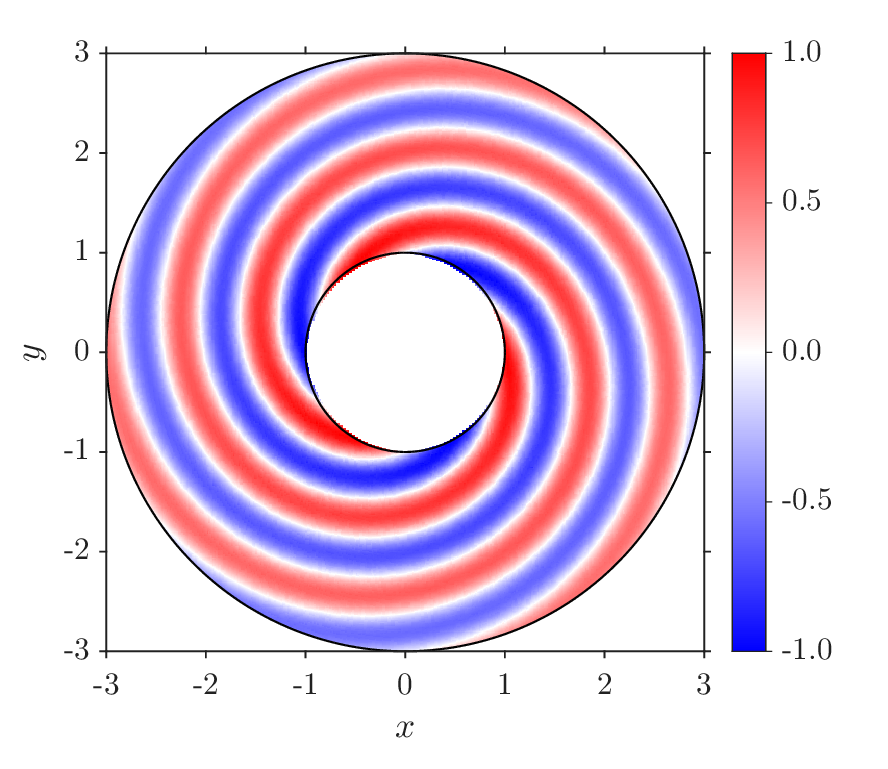}
\caption{Real part of the finite element solution, $\Re u_h$}
\label{fig:tbc_uhRe}
\end{figure}
\begin{figure}[H]\centering
\includegraphics[width=\mylenw]{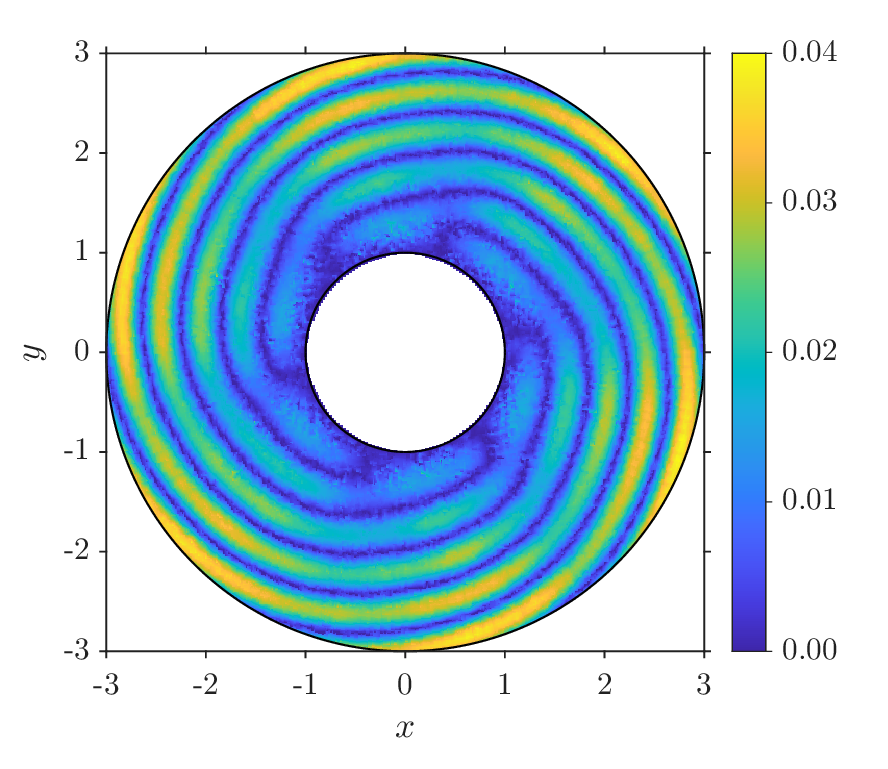}
\caption{Absolute error $|\Re u_h-\Re u|$}
\label{fig:tbc_errRe}
\end{figure}
In the computation, we set $\kappa=8$, $R_0=3$,
$N=500$, and $\rho=0.99R_0$.
Figure~\ref{fig:tbc_mesh} shows the finite element mesh used in the computation, which consists of 19,320 nodes and 37,974 linear triangular elements,
and the maximum mesh size is $h=0.0660$. Figure~\ref{fig:tbc_uexRe} displays the real part $\Re u$ of the exact solution, while Figure~\ref{fig:tbc_uhRe} shows the real part $\Re u_h$ of the finite element solution obtained using the proposed TBC.
The two plots are visually almost indistinguishable, indicating that the finite element solution successfully reproduces the main profile and oscillatory behavior of the exact solution.
This observation is supported by Figure~\ref{fig:tbc_errRe}, which shows that the absolute error
$|\Re u_h-\Re u|$
remains small over the whole computational domain, except for slight local increases near the artificial boundary. In particular, the maximum error is
\[
\|u_h-u\|_{C(\overline\Omega)}
\approx
\max_{1\leq j\leq M}
\left|u_h(\bm{x}_j)-u(\bm{x}_j)\right|= 0.0402, \]
where $M$ is the number of finite element nodes introduced in Section~\ref{sec:fem}, and
$\{\bm x_j\}_{j=1}^M\subset\overline\Omega$ denotes the set of finite element nodes.
The slight local increases near the artificial boundary remain small and do not indicate any significant deterioration in accuracy. This result suggests that the proposed TBC introduces no appreciable artificial reflections and accurately represents the exterior behavior.
These results confirm that the proposed TBC is correctly incorporated into the finite element formulation and provides an accurate and reliable boundary treatment for the exterior Helmholtz problem.
\subsection{Inhomogeneous problem with an irregular boundary}\label{sec:inhom}
\begin{figure}[H]\centering
\includegraphics[width=\mylenw]{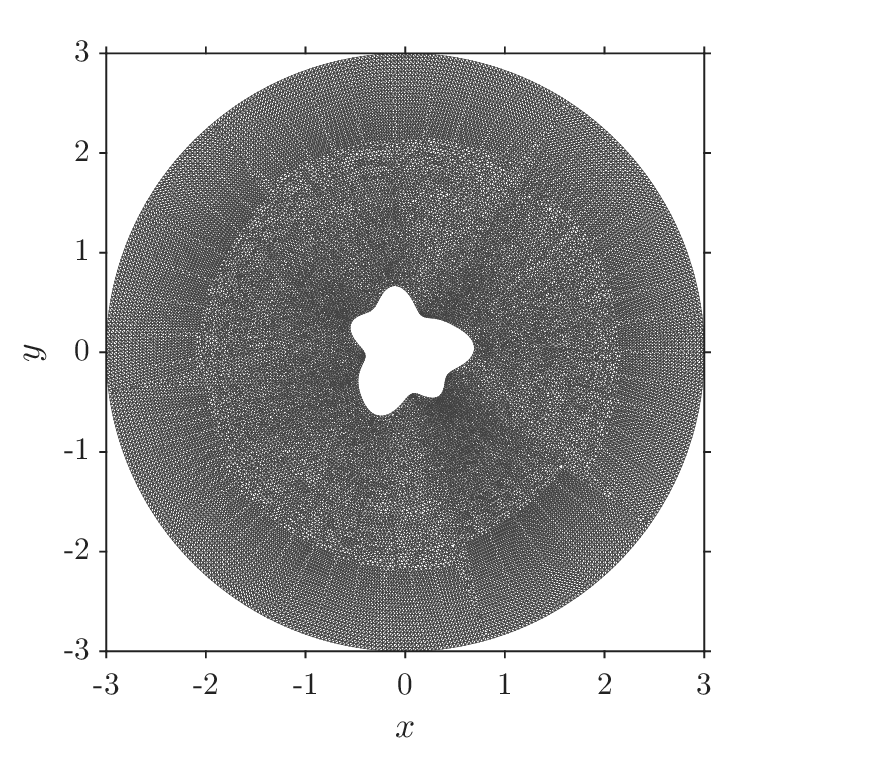}
\caption{Finite element mesh
(28,494 nodes and 56,322 elements)}
\label{fig:femesh}
\end{figure}

\begin{figure}[H]\centering
\includegraphics[width=\mylenw]{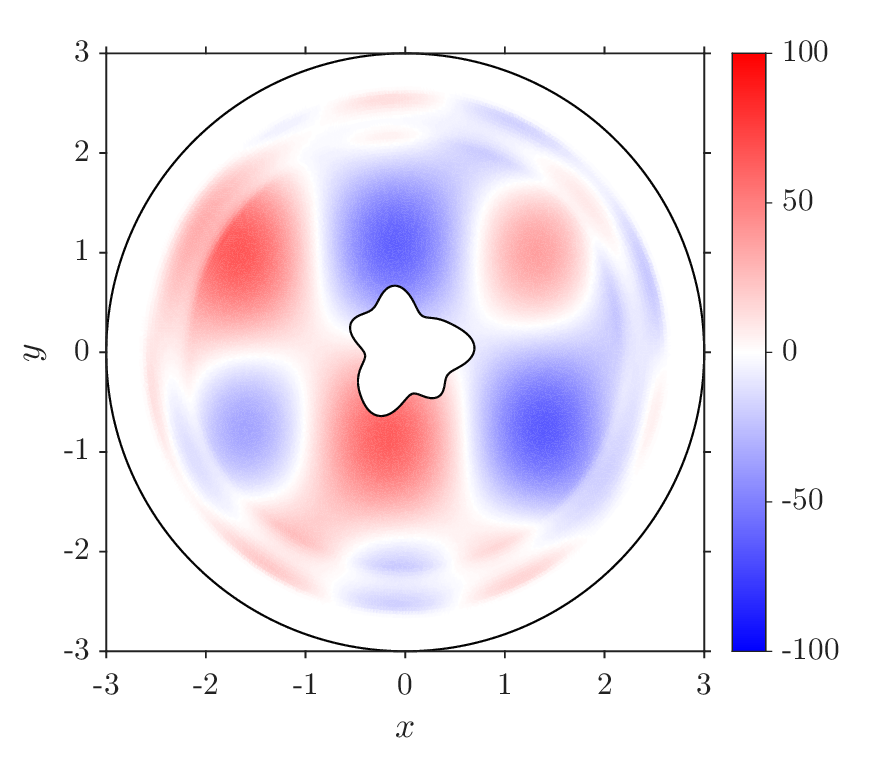}
\caption{Inhomogeneous term $f$}
\label{fig:inhomterm}
\end{figure}

\begin{figure}[H]\centering
\includegraphics[width=\mylenw]{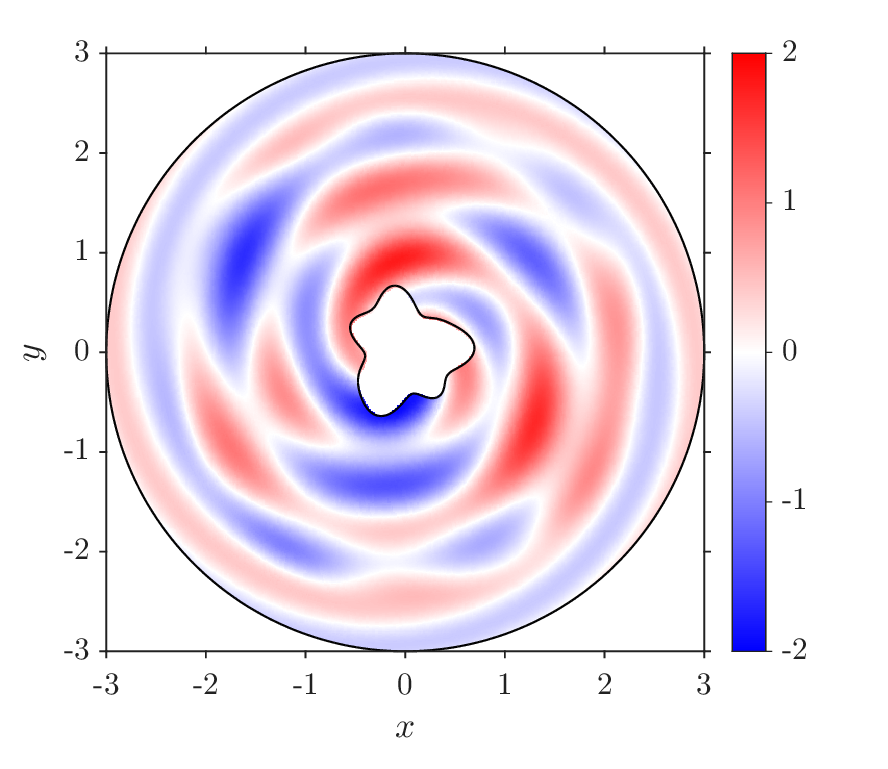}
\caption{Real part of the exact solution, $\Re u$}
\label{fig:fesol_ex}
\end{figure}

\begin{figure}[H]\centering
\includegraphics[width=\mylenw]{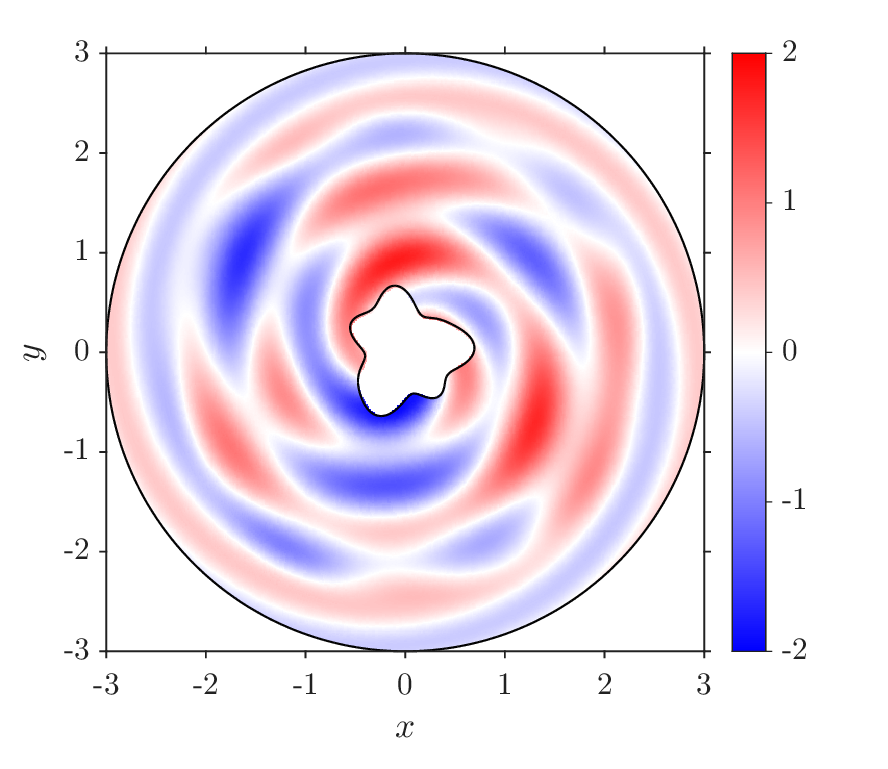}
\caption{Real part of the finite element solution, $\Re u_h$}
\label{fig:fesol}
\end{figure}

\begin{figure}[H]\centering
\includegraphics[width=\mylenw]{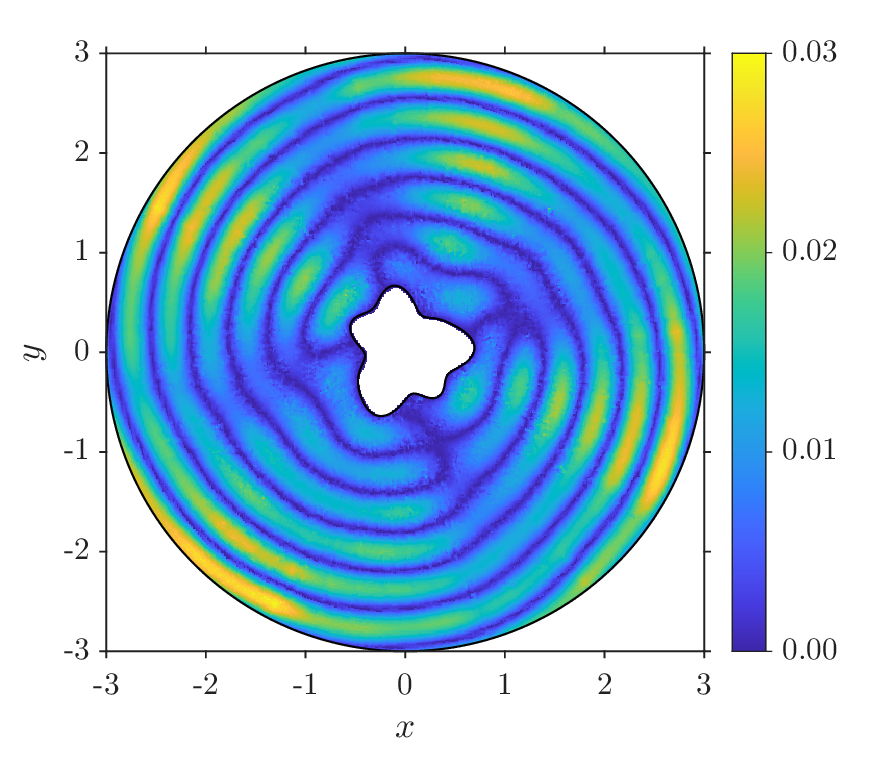}
\caption{Absolute error $|\Re u_h-\Re u|$}
\label{fig:feerr}
\end{figure}

This example is designed to verify that the finite element formulation incorporating the proposed MFS-based discrete DtN map remains accurate for an exterior Helmholtz problem with an irregular boundary and a compactly supported inhomogeneous term.
The boundary $\Gamma$ is given by
\[
R(\theta)
=
0.55
+0.10\cos(3\theta)
+0.06\sin(5\theta)
+0.04\cos(7\theta+0.3),
\qquad 0\le \theta < 2\pi .
\]
To validate the proposed method, we use a manufactured exact solution consisting of a compactly supported part and an outgoing radiating part:
\[
u=u_c+u_r,
\]
so that the corresponding inhomogeneous term
\[
f:=-(\Delta +\kappa^2) u
\]
is compactly supported in
the computational domain $\Omega$ enclosed by the obstacle boundary $\Gamma$ and the artificial boundary $\Gamma_0$ of radius $R_0=3$.
The detailed construction of $u$ and $f$ is given in  \ref{sec:apndx}.
In the computation, we set $\kappa=8$, $R_0=3$, $N=500$, and $\rho=0.99R_0$.
The resulting boundary value problem is discretized in $\Omega$ by the FEM using piecewise linear triangular elements.
Figure~\ref{fig:femesh} shows the finite element mesh used in the computation.
The resulting mesh consists of 28,494 nodes and 56,322  elements, and the maximum mesh size is $h=0.0671$.
The inhomogeneous term $f$, the real part $\Re u$ of the exact solution, the real part $\Re u_h$ of the finite element solution, and the absolute error $| \Re u_h-\Re u|$ are shown in Figures~\ref{fig:inhomterm}--\ref{fig:feerr}, respectively.
The finite element solution is in good agreement with the exact solution, indicating that the finite element approximation successfully captures the main features of the solution.
This observation is also supported by Figure~\ref{fig:feerr}, where the absolute error is seen to remain small throughout the computational domain, except for slight local increases near the artificial boundary and in regions where the solution varies rapidly.
Finally, the finite element solution $u_h$ is compared with the exact solution $u$ through the maximum error
\[
\|u_h-u\|_{C(\overline\Omega)}
\approx
\max_{1\leq j\leq M}
\left|u_h(\bm{x}_j)-u(\bm{x}_j)\right|=
0.0295 .
\]
This result confirms the effectiveness of the proposed method for exterior inhomogeneous Helmholtz problems with an irregular boundary.
\subsection{Acoustic propagation around a rigid barrier}
\label{sec:realprob}
\begin{figure}[H]
  \centering
  \begin{subfigure}[t]{0.495\linewidth}
    \centering
    \safeincludegraphics[width=\linewidth]{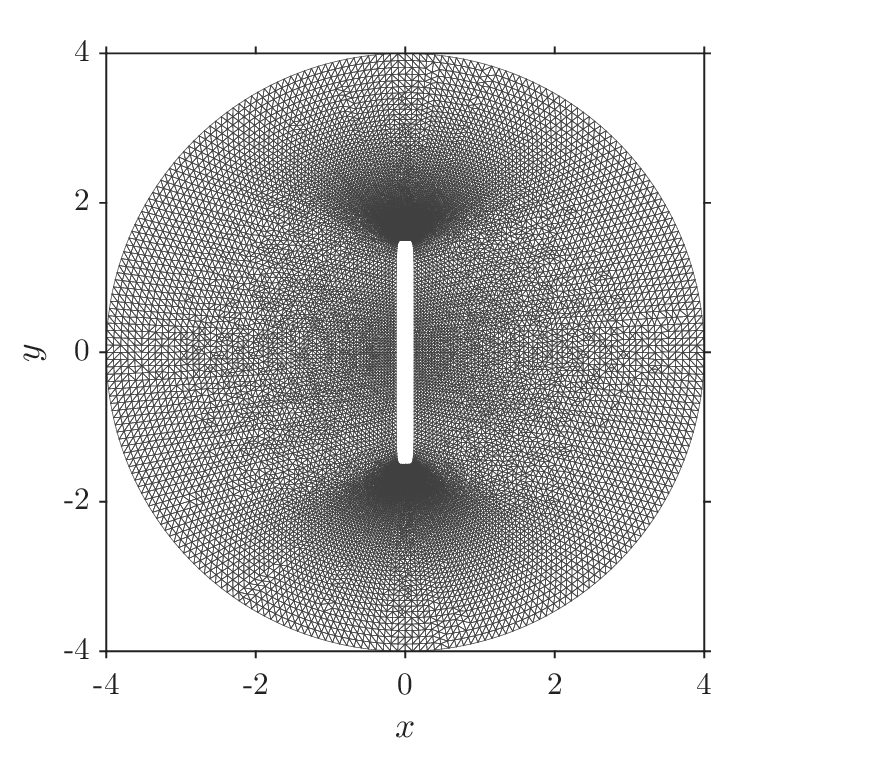}
    \caption{Finite element mesh (18,755 nodes and 36,994  elements)}
  \end{subfigure}
  \hfill
  \begin{subfigure}[t]{0.495\linewidth}
    \centering
    \safeincludegraphics[width=\linewidth]{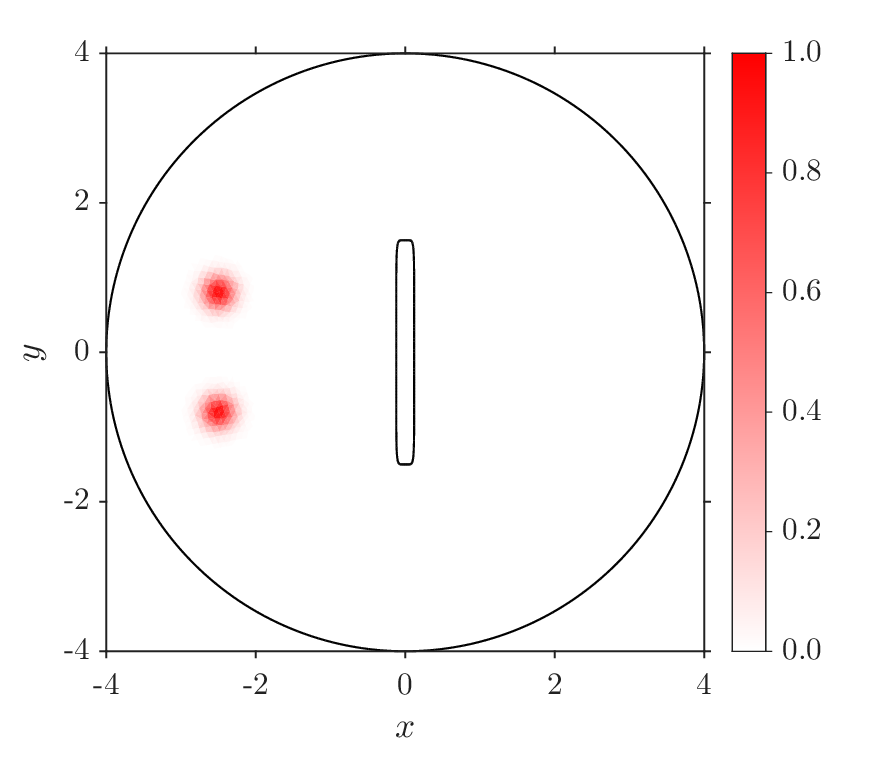}
    \caption{Inhomogeneous term $f$}
  \end{subfigure}
  \caption{Computational setting}
  \label{fig:mesh-source-en}
\end{figure}

\begin{figure}[H]
  \centering
  \begin{subfigure}[t]{0.495\linewidth}
    \centering
    \safeincludegraphics[width=\linewidth]{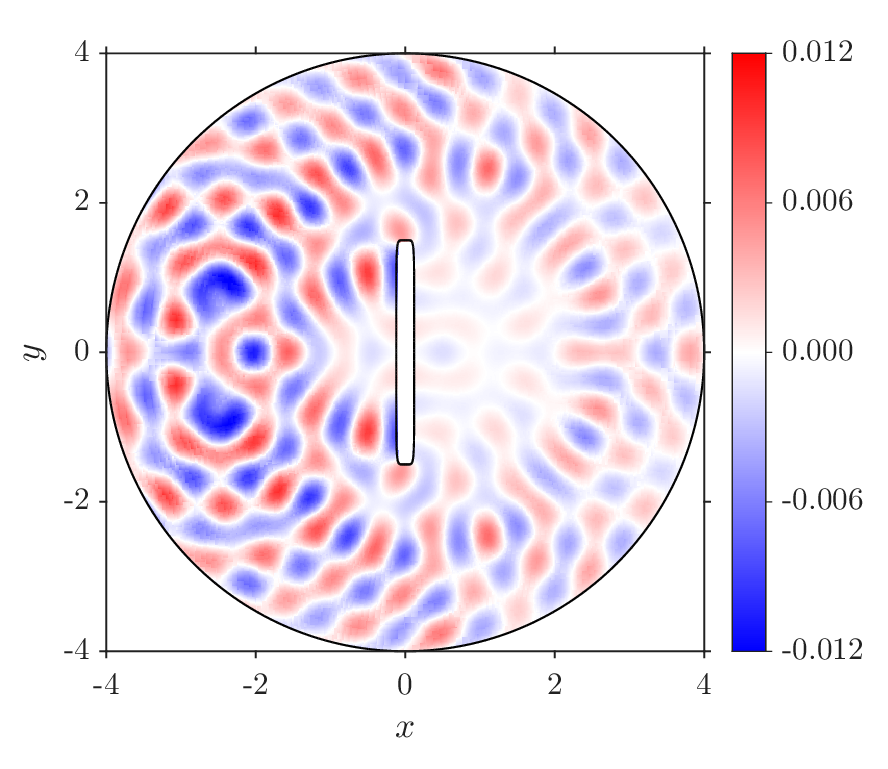}
    \caption{Real part of the finite element solution, $\Re u_h$}
  \end{subfigure}
  \hfill
  \begin{subfigure}[t]{0.495\linewidth}
    \centering
    \safeincludegraphics[width=\linewidth]{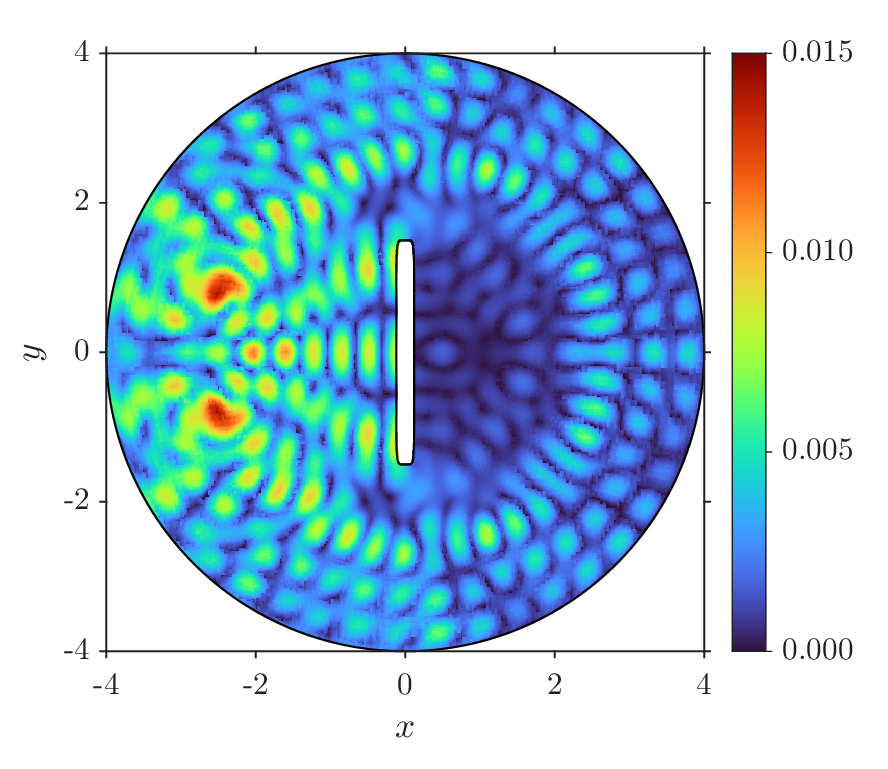}
    \caption{Absolute value of the finite element solution, $|u_h|$}
  \end{subfigure}
  \caption{Computed acoustic field around the rigid barrier}
  \label{fig:solution-fields-en}
\end{figure}

In this example, points in the plane are represented as vectors $\bm{x}\in\bR^2$, rather than as complex numbers $z\in\bC$.
We consider a two-dimensional acoustic propagation problem involving a finite rigid obstacle and two localized sound sources.
The obstacle is defined by
\[
\Omega_0
=
\left\{
(x,y)\in\bR^2
:
\left|\frac{x}{0.12}\right|^8
+
\left|\frac{y}{1.50}\right|^8
<1
\right\}.
\]
Then, $\Omega_0$ is a vertically elongated obstacle with rounded corners, measuring approximately $0.24\,\mathrm{m}$ in the horizontal direction and $3.0\,\mathrm{m}$ in the vertical direction.
The inhomogeneous term $f(\bm{x})$, representing the spatial distribution of the acoustic sources, is given by two in-phase Gaussian functions:
\[
f(\bm{x})
=
\exp\left(
-\frac{|\bm{x}-\bm{x}_{s,1}|^2}{2\sigma^2}
\right)
+
\exp\left(
-\frac{|\bm{x}-\bm{x}_{s,2}|^2}{2\sigma^2}
\right),
\]
where
\[
\bm{x}_{s,1}=(-2.5,0.8)\,\mathrm{m},
\qquad
\bm{x}_{s,2}=(-2.5,-0.8)\,\mathrm{m},
\qquad
\sigma=0.15\,\mathrm{m}.
\]
Although the Gaussian functions are not compactly supported in the strict mathematical sense, they are treated numerically as localized sources in the present computation.
The sound frequency and the sound speed are chosen as
$
\nu=500\,\mathrm{Hz}$ and
$c=343\,\mathrm{m/s}$,
so that
$\kappa={2\pi\nu}/{c}
\approx 9.16\,\mathrm{m}^{-1}$.
The corresponding wavelength is approximately \(0.686\,\mathrm{m}\).
The physical acoustic pressure associated with the complex amplitude \(u\) at position \(\bm{x}\) and time \(t\) is
\[
p(\bm{x},t)
=
\Re\left\{
u(\bm{x})e^{-2\pi i\nu t}
\right\}.
\]
A circular artificial boundary enclosing both the obstacle and the effective support of the inhomogeneous term is then introduced, with radius
$
R_0=4\,\mathrm{m}$.
In the theoretical formulation developed in the preceding sections, the Dirichlet condition is imposed on $\Gamma=\partial\Omega_0$.
In the present example, however, the obstacle is intended to represent an acoustically rigid barrier, and hence the homogeneous Neumann condition
\[
\frac{\partial u}{\partial n}=0
\quad\text{on}\quad\Gamma
\]
is imposed instead.
This modification does not affect the construction of the DtN map on the artificial boundary, because the DtN map depends only on the exterior Helmholtz problem outside the artificial boundary.
The difference lies only in how the boundary condition on $\Gamma$ is incorporated into the finite element formulation: the Dirichlet condition is incorporated into the finite element trial space, whereas the Neumann condition enters naturally through the boundary term in the weak formulation.
Figure~\ref{fig:mesh-source-en} shows the finite element mesh and the spatial distribution of the inhomogeneous term, whereas Figure~\ref{fig:solution-fields-en} shows the computed acoustic field.
The real part of the finite element solution exhibits the oscillatory structure generated by the two acoustic sources and the reflection from the rigid barrier.
Since the acoustic sources have equal amplitudes and phases and are placed symmetrically with respect to the \(x\)-axis, the resulting field is also approximately symmetric.
On the side of the barrier where the acoustic sources are located, the incident and reflected waves interfere, producing alternating regions of large and small amplitudes.
Behind the barrier, the direct wave is partially blocked, producing a region of
reduced amplitude.
The field does not vanish completely there because waves diffract around the upper and lower ends of the finite barrier.
The absolute value \(|u_h|\) clearly displays this reduction in amplitude behind the obstacle and is therefore useful for interpreting the shielding effect.
These results illustrate the reflection, interference, shielding, and edge diffraction produced by a finite rigid barrier in the presence of two in-phase sources.
\section{Conclusions}\label{sec:conc}
In this paper, we proposed a direct numerical method for solving the exterior Dirichlet problem for the two-dimensional inhomogeneous Helmholtz equation in an unbounded domain, based on a domain decomposition approach combined with a numerically constructed DtN map.
The DtN map is discretized by the MFS, leading to a circulant matrix representation that allows its efficient application via the FFT. By incorporating the TBC defined by the DtN map into the finite element discretization, the original problem in the unbounded domain is reduced to a problem posed in a bounded computational domain. In this way, the proposed method enables the FEM to handle exterior problems that cannot be treated directly by a standard finite element formulation in a bounded domain alone. Furthermore, since the present framework also applies to inhomogeneous equations, it overcomes a difficulty of the classical boundary element method, whose application becomes less straightforward in the presence of inhomogeneous terms.
The numerical examples serve distinct purposes: they demonstrate that the MFS-based discrete DtN map can be computed accurately even for a large wavenumber, confirm the computational efficiency of the FFT-based algorithm for large numbers of collocation points, validate the finite element formulation incorporating the discrete DtN map, and illustrate the applicability of the proposed method to both a manufactured inhomogeneous problem and a physically motivated one.
From a theoretical viewpoint,
it is known that the DtN map defined on a circular boundary
is diagonal with respect to the Fourier basis \cite{BTL},
that is, the Fourier basis forms an eigenfunction system of the DtN map.
This property corresponds to the fact that the circulant matrix
representing the discrete DtN map obtained by the MFS
is diagonalized by the DFT matrix.
From this perspective,
it is also of interest to investigate in what topology
the matrix representing the discrete DtN map converges
to the exact DtN map
as the number of collocation points tends to infinity.
Furthermore, when three-dimensional problems are considered,
the matrices arising in the MFS discretization are no longer simple circulant matrices.
Instead, they have a block structure
whose components are circulant matrices
\cite{Liu-c2}.
Therefore, the idea of the proposed method
can be extended to three-dimensional problems.
Related MDAs for the MFS on spherical geometries are described in Section~9.2.2 of \cite{CCK}, and a numerical example for the three-dimensional Helmholtz equation  is presented in Example~9.4 of the same reference. Nevertheless, the efficient explicit construction of the discrete DtN map for incorporation into a finite element formulation remains a nontrivial computational issue in three dimensions.
More specifically, the resulting matrices are considerably larger, and constructing an explicit matrix representation of the discrete DtN map generally requires solving large linear systems with multiple right-hand sides.
Equivalently, it requires evaluating the action of the inverse of a large matrix on multiple vectors.
Therefore, further computational techniques are required to make the three-dimensional extension practically efficient. This issue will be investigated in future work.
\appendix
\section{Construction of the manufactured solution for the inhomogeneous problem}
\label{sec:apndx}
We describe the construction of the manufactured solution
used in Section~\ref{sec:inhom} for the inhomogeneous exterior Helmholtz problem
with an irregular boundary.
Although the main text identifies $\mathbb{R}^2$ with $\mathbb{C}$, in this appendix we describe the manufactured solution in Cartesian coordinates $(x,y)\in\mathbb{R}^2$ for clarity.
We write the solution in the form
\[
u=u_c+u_r,
\]
where $u_c$ is introduced to generate a compactly supported inhomogeneous term,
whereas $u_r$ is an outgoing homogeneous solution ensuring that the solution
does not vanish outside the artificial boundary.
\subsection*{Compactly supported part}
The first part is defined by
\[
u_c(x,y)
:=
\chi(r)\,p(x,y),
\qquad
r:=\sqrt{x^2+y^2},
\]
where $p$ is the smooth real-valued function
\[
p(x,y)
:=
\cos(2.1x+0.3)\sin(1.7y-0.2)+0.15x+0.10y.
\]
The cutoff function $\chi$ is chosen so that it is equal to $1$ near the boundary $\Gamma$
and vanishes near the artificial boundary $\Gamma_0$.
More precisely, for two radii $0<R_1<R_2<R_0$, we define
\[
\chi(r)
:=
\begin{cases}
1, & 0\le r\le R_1,\\[1mm]
1-S\!\left(\dfrac{r-R_1}{R_2-R_1}\right), & R_1<r<R_2,\\[3mm]
0, & r\ge R_2,
\end{cases}
\]
where
\[
S(s):=10s^3-15s^4+6s^5.
\]
In the computation, we take
\[
R_1=0.72R_0,
\qquad
R_2=0.88R_0.
\]
Since $\chi(r)=0$ for $r\ge R_2$, the function $u_c$ vanishes
in a neighborhood of the artificial boundary $\Gamma_0$.
\subsection*{Outgoing radiating part}
The second part is chosen as a single outgoing Fourier--Hankel mode:
\[
u_r(r\cos\theta,r\sin\theta)
:=
\beta\,
\frac{H_m^{(1)}(\kappa r)}{H_m^{(1)}(\kappa R_0)}
e^{im\theta},
\]
where, as stated in the main text, $H_m^{(1)}$ denotes the Hankel function of the first kind of order $m$.
In the present example, the parameters are chosen as
\[
m=2,
\qquad
\beta=0.35+0.20\,i.
\]
By construction, $u_r$ satisfies
\[
(\Delta +\kappa^2) u_r=0
\]
in the computational domain and represents an outgoing wave outside the artificial boundary.
Our aim is to construct $u$ so that the resulting
inhomogeneous term is compactly supported, while $u$ itself does not vanish outside $\Gamma_0$.
\subsection*{Inhomogeneous term}
The inhomogeneous term is defined by
\[
f:=-(\Delta +\kappa^2)u.
\]
Since $u_r$ satisfies the homogeneous Helmholtz equation, we obtain
\[
f=-(\Delta +\kappa^2)u_c.
\]
Substituting $u_c=\chi p$, we have
\[
f
=-\left(
\Delta+\kappa^2
\right)(\chi p)
=-\left(
p\,\Delta\chi
+2\nabla\chi\cdot\nabla p
+\chi\,\Delta p
+\kappa^2\chi p\right).
\]
The derivatives of $p$ are given by
\[
\nabla p
=
\begin{pmatrix}
-2.1\sin(2.1x+0.3)\sin(1.7y-0.2)+0.15\\[1mm]
1.7\cos(2.1x+0.3)\cos(1.7y-0.2)+0.10
\end{pmatrix},
\]
and
\[
\Delta p
=
-(2.1^2+1.7^2)\cos(2.1x+0.3)\sin(1.7y-0.2).
\]
Since $\chi$ depends only on $r=\sqrt{x^2+y^2}$, we also have
\[
\nabla\chi(r)
=
\chi'(r)\frac{1}{r}
\begin{pmatrix}
x\\ y
\end{pmatrix},
\qquad
\Delta\chi(r)
=
\chi''(r)+\frac{1}{r}\chi'(r),
\qquad r>0.
\]
Let $s:=(r-R_1)/(R_2-R_1)$.
Since
\[
S'(s)=30s^2-60s^3+30s^4,\qquad
S''(s)=60s-180s^2+120s^3,
\]
the first and second derivatives of the cutoff function are explicitly
given by
\[
\chi'(r)
=
\begin{cases}
-\dfrac{30s^2-60s^3+30s^4}{R_2-R_1},
& R_1<r<R_2,\\[3mm]
0,
& 0\leq r\leq R_1 \text{ or } r\geq R_2,
\end{cases}
\]
and
\[
\chi''(r)
=
\begin{cases}
-\dfrac{60s-180s^2+120s^3}{(R_2-R_1)^2},
& R_1<r<R_2,\\[3mm]
0,
& 0\leq r\leq R_1 \text{ or } r\geq R_2.
\end{cases}
\]
Consequently, for $r>0$, the inhomogeneous term can be evaluated explicitly as
\[
\begin{aligned}
f(x,y)
={}&
-p(x,y)
\left(
\chi''(r)+\frac{\chi'(r)}{r}
\right)\\
&-\frac{2\chi'(r)}{r}
\left(
x\frac{\partial p}{\partial x}(x,y)
+
y\frac{\partial p}{\partial y}(x,y)
\right)\\
&-\chi(r)\Delta p(x,y)
-\kappa^2\chi(r)p(x,y).
\end{aligned}
\]
Thus, all the quantities on the right-hand side are explicitly specified,
and the inhomogeneous term used in the computation is uniquely determined.
Hence, the support of $f$ is contained in the closed circular region\[
\operatorname{supp}f
\subset
\{(x,y)\in\mathbb{R}^2: \sqrt{x^2+y^2}\le R_2\}, \]
so that $f$ is compactly supported in the computational domain.
\subsection*{Boundary data}
The Dirichlet boundary data on $\Gamma$ are taken from the solution,
\[
g=u|_{\Gamma}.
\]
On the artificial boundary $\Gamma_0$, the TBC is imposed.
Since $u_c$ vanishes near $\Gamma_0$, both $u$ and $\partial u/\partial n$ on $\Gamma_0$ are determined solely by $u_r$.
\section*{Acknowledgments}
The author would like to express his sincere gratitude to the reviewers for their careful reading and constructive and helpful comments, which greatly improved the quality of the manuscript.

\end{document}